\def \NN {\mathbb N}
\def \CC {\mathbb C}

\def \RR {\mathbb R}
\def \ZZ {\mathbb Z}
\def \HH {\mathbb H}

\def \epsilon{\varepsilon}

\def \A  {{\mathcal A}}

\def \C  {{\mathcal C}}
\def \D  {{\mathcal D}}
\def \E  {{\mathcal E}}

\def \LL {{\mathcal L}}

\def \S  {{\mathcal S}}

\def \d {\text{d}}

\def \bfx {{\boldsymbol{x}}}

\def \bfw {{\boldsymbol{w}}}

\def \half {{\textstyle{ 1\over 2}}}

\def \ep {\epsilon}
\def \ga {\gamma}
\def \Ga {\Gamma}
\def \si {\sigma}
\def \la {\lambda}

\newcommand{\res}{{\rm res}}
\renewcommand{\S}{{\mathcal S}}

\documentclass[12pt,reqno]{amsart}

\usepackage{amsmath,amssymb}

\usepackage{url} 

\numberwithin{equation}{section}

\begin{document}


\title[]{Classification of $L$-functions of degree $2$ \\ and conductor $1$} 

\author[]{J.KACZOROWSKI \lowercase{and} A.PERELLI}
\maketitle

{\bf Abstract.} We give a full description of the functions $F$ of degree 2 and conductor 1 in the general framework of the extended Selberg class $\S^\sharp$. This is performed by means of a new numerical invariant $\chi_F$, which is easily computed from the data of the functional equation. We show that the value of $\chi_F$ gives a precise description of the nature of $F$, thus providing a sharp form of the classical converse theorems of Hecke and Maass. In particular, our result confirms, in the special case under consideration, the conjecture that the functions in the Selberg class $\S$ are automorphic $L$-functions.



\bigskip
\centerline{Contents}

1. Introduction

2. Definitions, notation and basic requisites

3. Outline of the proof

4. Invariants

5. Virtual $\ga$-factors

6. Period functions

7. Conclusion of the proof

References

\vskip.5cm
\section{Introduction}

\smallskip
It is generally expected that the Selberg class $\S$ of $L$-functions coincides with the class of automorphic $L$-functions. This is known for functions of degrees $0<d<2$, see Conrey-Ghosh \cite{Co-Gh/1993} and the authors' papers \cite{Ka-Pe/1999a} and \cite{Ka-Pe/2011}, but already the important case of degree $d=2$ appears at present to be out of reach in its full generality. In that case, it is conjectured that the primitive functions in $\S$ coincide, roughly, with the $L$-functions of holomorphic and non-holomorphic eigenforms of any level.  Since the Ramanujan conjecture is not known in the non-holomorphic case, while it is assumed in the definition of $\S$, here we consider the weaker conjecture that the functions of degree 2 in $\S$ are contained, roughly, in the family of the $L$-functions of the above eigenforms. We refer to Section 2 for the definitions and notions introduced in this section. 

\smallskip
If the functions $F\in\S$ satisfy a functional equation with very special $\Ga$-factors and, if the conductor $q$ is greater than 4, certain twists of $F$ by Dirichlet characters also satisfy suitable functional equations, then the latter conjecture is known thanks to the classical converse theorems of Hecke \cite{Hec/1983}, Maass \cite{Maa/1949}, Weil \cite{Wei/1967}, Jacquet-Langlands \cite{Ja-La/1970} and the theory of Hecke operators. Concerning general $\Ga$-factors, the only known case of such a conjecture is, as far as we know, when the conductor $q$ equals 1 and $F$ has a pole at $s=1$, in which case it turns out that necessarily $F=\zeta^2$; see \cite{Ka-Pe/2015}. Although the existence of an Euler product expansion for $F$ is not necessary in the converse theorems of Hecke, Maass and Weil, both the Euler product and the pole at $s=1$ play a crucial role in the arguments in \cite{Ka-Pe/2015}.

\smallskip
In this paper we describe the nature of the functions of degree $d=2$ and conductor $q=1$ in the full generality of the extended Selberg class $\S^\sharp$ and, as a consequence, we confirm the above conjecture in the case of conductor 1. Our description of such functions is performed by means of a new numerical invariant, denoted $\chi_F$ and called the {\it eigenweight} of $F\in\S^\sharp$, which is easy to compute from the data of the functional equation. The name of $\chi_F$ comes from the fact that it determines the eigenvalue and the weight of the modular form corresponding to $F$; see Theorem 1.1. To give a glimpse on our main result, let $F\in\S^\sharp$ with $d=2$ and $q=1$ be normalized in the sense described below; then, for example,
\[
\chi_F=0 \ \Longrightarrow \ F(s)= \zeta(s)^2 \quad \text{and} \quad  \chi_F = \frac{121}{2} \Longrightarrow F(s) = L(s +\frac{11}{2},\Delta),
\]
where $L(s,\Delta)$ is the Ramanujan $L$-function.

\smallskip
In order to state our results we have to introduce a normalization. This is due to the fact that the functional equation in $\S^\sharp$ reflects $F$ into its conjugate $\overline{F}$, while the $L$-functions of modular forms of level 1 are reflected into themselves. We say that $F\in\S^\sharp$ is {\it normalized} if its internal shift $\theta_F$ vanishes and the first nonvanishing Dirichlet coefficient is 1. Normalized functions have a twofold advantage, making them the right objects for our purpose. Indeed, on one side every $F\in\S^\sharp$ with $d=2$ and $q=1$ can be normalized by means of a simple procedure, see Lemma 4.1, so we may consider only such functions without loosing generality. On the other side, normalized functions have real coefficients, see again Lemma 4.1, and hence the above mentioned discrepancy disappears.

\smallskip
The new invariant $\chi_F$ is defined as
\begin{equation}
\label{1-1}
\chi_F = \xi_F+H_F(2) +2/3,
\end{equation}
where $ \xi_F$ and $H_F(2)$ are the $H$-invariants defined in \eqref{2-3}-\eqref{2-5} below. Hence $\chi_F$ can easily be computed plugging the data of the functional equation of $F$ into the first and second Bernoulli polynomials. For example, if  $f$ is a holomorphic cusp form of level 1 and weight $k$, with first nonvanishing Fourier coefficient equal 1, and we denote by $L(s,f)$ its $L$-function, then $F(s)=L(s+\frac{k-1}{2},f)$ is a normalized function of $\S^\sharp$ with $d=2$ and $q=1$, and
\[
\chi_F = \frac{(k-1)^2}{2}.
\]
Similarly, one easily checks that if $u$ is a Maass form of level 1 and weight 0, with eigenvalue $1/4+\kappa^2$ and first Fourier-Bessel coefficient equal 1, and $L(s,u)$ is its $L$-function, then $F(s)=L(s,u)$ is a normalized member of $\S^\sharp$ with $d=2$ and $q=1$, and
\[
\chi_F = - 2\kappa^2.
\]

\smallskip
Our main result shows, conversely, that the value of $\chi_F$ is sufficient to detect the nature of any normalized $F\in\S^\sharp$ of degree 2 and conductor 1. Thus it may be regarded as a sharp form of the Hecke and Maass converse theorems; note that these theorems are used in the very last step of our proof.

\medskip
{\bf Theorem 1.1.} {\sl Let $F\in\S^\sharp$ of degree $2$ and conductor $1$ be normalized. Then $\chi_F\in\RR$ and

\smallskip
(i) if $\chi_F>0$ then there exists a holomorphic cusp form $f$ of level $1$ and even integral weight $k= 1 + \sqrt{2\chi_F}$ such that $F(s) = L(s+\frac{k-1}{2},f)$;

\smallskip
(ii) if $\chi_F=0$ then $F(s) = \zeta(s)^2$;

\smallskip
(iii) if $\chi_F<0$ then there exists a Maass form $u$ of level $1$, weight $0$ and with eigenvalue $1/4+\kappa^2= (1-2\chi_F )/4$ such that $F(s) = L(s,u)$.}

\medskip
In case (iii) we can specify the parity $\ep$ of $u$ by means of the root number $\omega_F$ of $F$, namely 
\begin{equation}
\label{1-2}
\ep= \frac{1-\omega_F}{2}.
\end{equation}
Clearly, if $F\in\S^\sharp$ with $d=2$ and $q=1$ is not normalized, we may first normalize it and then use Theorem 1.1 to detect its nature. Therefore every such $F$ is closely related to one of the three types of $L$-functions in Theorem 1.1.

\medskip
{\bf Corollary 1.1.} {\sl Every $F\in\S^\sharp$ of degree $2$ and conductor $1$ belongs, modulo normalization, to one of the three families of $L$-functions in Theorem $1.1$.}

\medskip
Concerning the functions in $\S$, from the theory of Hecke operators and the linear independence of Euler products in Kaczorowski-Molteni-Perelli \cite{K-M-P/1999},\cite{K-M-P/2006} we have the following corollary.

\medskip
{\bf Corollary 1.2.}  {\sl Every $F\in\S$ of degree $2$ and conductor $1$ is an automorphic $L$-function.}

\medskip
The methods in \cite{Ka-Pe/1999a},\cite{Ka-Pe/2011} and \cite{Ka-Pe/2015} are apparently not sufficient to prove Theorem 1.1, and genuinely new ideas needed to be introduced. We may briefly summarise such ideas as follows; we always assume that $F\in\S^\sharp$ has degree 2, conductor 1 and is normalized.

(a) First, essential use is made of the {\it structural invariants} $d_\ell(F)$, see \eqref{2-11}, the crucial new result being that for $N\geq2$ the points $\big(d_1(F),\dots,d_N(F)\big)$ lie on certain universal algebraic varieties. We suspect that something similar should hold in general for the functions $F\in\S^\sharp$, i.e. the $d_\ell(F)$ should lie on certain algebraic varieties to a large extent independent of $F$. In turn, this could explain why the $L$-functions satisfy only functional equations with very special $\Ga$-factors, and in particular could shed some light on the general structure of the Selberg class.

(b) Thanks to (a), we can associate to any function $F$ a unique {\it virtual $\ga$-factor} of classical shape, see \eqref{5-1},  and show that $F$ satisfies a functional equation of a slightly different type, but involving such a virtual $\ga$-factor.

(c) Finally we show that the original $\ga$-factor of $F$ actually coincides with the above virtual $\ga$-factor; this is performed studying certain {\it period functions}, in the sense of Lewis-Zagier \cite{Le-Za/2001}, associated with $F$; see \eqref{6-1}. 

Now the theorem follows at once from the classical converse theorems by Hecke and Maass. We refer to Section 3 for a more detailed outline of the proof, and conclude this section with two remarks.

\medskip
{\bf Remark 1.1.} Since it is well known that exceptional eigenvalues do not exist in level 1, see Selberg \cite{Sel/1965}, we use this fact in the proof. However, if exceptional eigenvalues would actually exist, then our proof would still work, and the resulting theorem would list a second possibility for the functions $F$ with $\chi_F<0$, namely $F(s)=L(s,u)$ for some Maass form $u$ associated with an exceptional eigenvalue. \qed

\medskip
{\bf Remark 1.2.} Long ago we conjectured that the functional equation of any $F\in\S^\sharp$ of degree $d$ is completely described by conductor, root number and the $H$-invariants $H_F(n)$ with $n\leq d$, see \cite[p.103]{Ka-Pe/2002}. This is confirmed by Theorem 1.1 and \eqref{1-2} in the present special case where $d=2$ and $q=1$, since $\chi_F$ is defined in terms of $H_F(1)$ and $H_F(2)$; see \eqref{1-1} and \eqref{2-5}. \qed

\medskip
{\bf Acknowledgements.} This research was partially supported by the Istituto Nazionale di Alta Matematica, by the MIUR grant PRIN-2017 {\sl ``Geometric, algebraic and analytic methods in arithmetic''} and by grant 2017/25/B/ST1/00208 {\sl ``Analytic methods in number theory''}  from the National Science Centre, Poland.

\medskip
\section{Definitions, notation and basic requisites}

\smallskip
{\bf 2.1. Definitions and notation.}
Throughout the paper we write $s=\si+it$ and $\overline{f}(s)$ for $\overline{f(\overline{s})}$, $f(s)\equiv0$ means that $f(s)$ vanishes identically and $\HH$ denotes the upper half-plane $\{z=x+iy\in\CC:y>0\}$. The extended Selberg class $\S^\sharp$ consists of non identically vanishing Dirichlet series 
\[
F(s)= \sum_{n=1}^\infty \frac{a(n)}{n^s},
\]
absolutely convergent for $\si>1$, such that $(s-1)^mF(s)$ is entire of finite order for some integer $m\geq0$, and satisfying a functional equation of type
\begin{equation}
\label{2-1}
F(s) \gamma(s) = \omega \overline{\gamma}(1-s) \overline{F}(1-s),
\end{equation}
where $|\omega|=1$ and the $\gamma$-factor
\begin{equation}
\label{2-2}
\gamma(s) = Q^s\prod_{j=1}^r\Gamma(\lambda_js+\mu_j) 
\end{equation}
has $Q>0$, $r\geq0$, $\lambda_j>0$ and $\Re(\mu_j)\geq0$. Note that the conjugate function $\overline{F}$ has conjugated coefficients $\overline{a(n)}$. The Selberg class $\S$ is, roughly speaking, the subclass of $\S^\sharp$ of the functions with Euler product and satisfying the Ramanujan conjecture $a(n)\ll n^\ep$. We refer to Selberg \cite{Sel/1989}, Conrey-Ghosh \cite{Co-Gh/1993} and to our survey papers \cite{Kac/2006},\cite{Ka-Pe/1999b},\cite{Per/2005},\cite{Per/2004},\cite{Per/2010},\cite{Per/2017} for further definitions, examples and the basic theory of the Selberg class. 

\smallskip
Degree $d$, conductor $q$, root number $\omega_F$ and $\xi$-invariant $\xi_F$ of $F\in\S^\sharp$ are defined by
\begin{equation}
\label{2-3}
\begin{split}
d=d_F :=2\sum_{j=1}^r\lambda_j, \qquad q=q_F:= (2\pi)^dQ^2\prod_{j=1}^r\lambda_j^{2\lambda_j}, \\
\omega_F=\omega \prod_{j=1}^r \lambda_j^{-2i\Im(\mu_j)},  \qquad \xi_F = 2\sum_{j=1}^r(\mu_j-1/2):= \eta_F+ id\theta_F
\end{split}
\end{equation}
with $\eta_F,\theta_F\in\RR$; $\theta_F$ is called the internal shift, and the classical $L$-functions have $\theta_F=0$. In this paper we deal only with $F\in\S^\sharp$ of positive degree, hence with $r\geq 1$. The $H$-invariants, introduced in \cite{Ka-Pe/2002}, are defined for every $n\geq0$ as
\begin{equation}
\label{2-4}
H_F(n) = 2 \sum_{j=1}^r \frac{B_n(\mu_j)}{\lambda_j^{n-1}},
\end{equation}
where $B_n(z)$ is the $n$th Bernoulli polynomial. Clearly
\begin{equation}
\label{2-5}
H_F(0)=d, \qquad H_F(1)=\xi_F,  \qquad H_F(2) = 2 \sum_{j=1}^r \frac{\mu_j^2 -\mu_j + 1/6}{\lambda_j}.
\end{equation}

\smallskip
{\bf Remark 2.1.}  Since $H$-invariants depend only on the data $Q, r, \lambda_j,\mu_j$ of $\gamma$-factors \eqref{2-2}, we may define such invariants for any $\gamma$-factor, without referring to functions $F\in\S^\sharp$. More generally, the same remark holds for any invariant depending only on the data of $\gamma$-factors. Clearly, the invariants of a $\gamma$-factor coincide with the corresponding invariants of any $F\in\S^\sharp$ having such a $\gamma$-factor. The invariants of $\gamma$-factors are usually denoted replacing the suffix $F$ by $\gamma$. \qed

\smallskip
For $\alpha>0$ the standard twist of $F\in\S^\sharp$ is
\[
F(s,\alpha) = \sum_{n=1}^\infty \frac{a(n)}{n^s} e(-\alpha n^{1/d}), \quad e(x) = e^{2\pi i x},
\]
the spectrum of $F$ is defined as
\begin{equation}
\notag
\text{Spec}(F) = \{\alpha>0: a(n_\alpha)\neq0\} = \Big\{d\big(\frac{m}{q}\big)^{1/d}: m\in\NN \ \text{with} \ a(m)\neq 0\Big\},
\end{equation}
where
\begin{equation}
\notag
n_\alpha = q d^{-d} \alpha^{d} \quad \text{and} \quad a(n_\alpha)=0 \ \text{if} \ n_\alpha\not\in \NN.
\end{equation}
Moreover, for $\ell=0,1,\dots$ we write
\begin{equation}
\label{2-6}
s_\ell = \frac{d+1}{2d} -\frac{\ell}{d} \quad \text{and} \quad s^*_\ell = s_\ell -i\theta_F. 
\end{equation}
Finally we consider the $S$-function
\begin{equation}
\label{2-7}
S_F(s) := 2^r \prod_{j=1}^r \sin(\pi(\lambda_js + \mu_j)) = \sum_{j=-N}^N a_j e^{i\pi d_F\omega_j s}
\end{equation}
with certain $N\in\NN$, $a_j\in\CC$ and $-1/2 =\omega_{-N} < \dots <\omega_N = 1/2$, and the $h$-function
\[
h_F(s) = \frac{\omega}{(2\pi)^r} Q^{1-2s} \prod_{j=1}^r \big(\Gamma(\lambda_j(1-s)+\overline{\mu}_j) \Gamma(1-\lambda_js-\mu_j)\big).
\]
Here we keep the notation for $S_F(s)$ used in \cite{Ka-Pe/revisited} since we heavily refer to that paper in Section 4, but in Section 5 we make a slight change of notation in \eqref{2-7}, see \eqref{5-14}. Since the function $h_F(s)$ depends also on the $\omega$-datum, not present in $\gamma$-factors, we define the analog of $h_F(s)$ for $\gamma$-factors as
\begin{equation}
\label{2-8}
h_\gamma(s) =  \prod_{j=1}^r \lambda_j^{2i\Im(\mu_j)} \frac{1}{(2\pi)^r} Q^{1-2s} \prod_{j=1}^r \big(\Gamma(\lambda_j(1-s)+\overline{\mu}_j) \Gamma(1-\lambda_js-\mu_j)\big).
\end{equation}
Instead, the analog $S_\gamma(s)$ of $S_F(s)$ is defined exactly in the same way. Clearly, in view of \eqref{2-3}, for any $\gamma$-factor of $F$ we have that
\begin{equation}
\label{2-9}
h_F(s)= \omega_F h_\gamma(s).
\end{equation}

\medskip
{\bf 2.2. Basic requisites.}
We recall that every $F\in\S^\sharp$ has polynomial growth on vertical strips. Moreover, the standard twist $F(s,\alpha)$ is entire if $\alpha\not\in$ Spec$(F)$, while for $\alpha\in$ Spec$(F)$ it is meromorphic on $\CC$ with at most simple poles at the points $s=s^*_\ell$, with residue denoted by $\rho_\ell(\alpha)$. It is known that $\rho_0(\alpha)\neq0$ for every $\alpha\in$ Spec$(F)$. Further, $F(s,\alpha)$ has polynomial growth on every vertical strip. We refer to \cite{Ka-Pe/2005} and \cite{Ka-Pe/2016a} for these and other results on $F(s,\alpha)$.

\smallskip
The functions $S_F(s)$, $S_\gamma(s)$, $h_F(s)$ and $h_\gamma(s)$ are invariants, and functional equation \eqref{2-1} can be written in the invariant form
\begin{equation}
\label{2-10}
F(s) = h_F(s) S_F(s) \overline{F}(1-s),
\end{equation}
see Section 1.2 of \cite{Ka-Pe/revisited}. Moreover, for large $|s|$ outside an arbitrarily small angular region containing the positive real axis we have the asymptotic expansion
\begin{equation}
\label{2-11}
h_F(s)  \approx \frac{\omega_F}{\sqrt{2\pi}} \left(\frac{q^{1/d}}{2\pi d}\right)^{d(\frac12-s)} \sum_{\ell=0}^\infty d_\ell(F) \Gamma\big(d(s_\ell^*-s)\big), \qquad d_0(F)= d^{id\theta_F},
\end{equation}
where $ \approx$ means that cutting the sum at $\ell=M$ one gets a meromorphic remainder which is $\ll$ than the modulus of the $M$-th term times $1/|s|$; see Section 1.3 of \cite{Ka-Pe/revisited}. Since \eqref{2-11} is proved by means of Stirling's formula without using the properties of $F\in\S^\sharp$, see Section 3.2 of \cite{Ka-Pe/revisited}, a completely analogous expansion holds for $h_\gamma(s)$ as well. Clearly, such an expansion holds without the factor $\omega_F$ and with $d_\ell(F)$ replaced by $d_\ell(\gamma)$. The coefficients $d_\ell(F)$ and $d_\ell(\gamma)$ are called the {\it structural invariants} and play an important role in the Selberg class theory. When $\alpha\in$ Spec$(F)$, Theorem 3 of \cite{Ka-Pe/revisited} shows that the residue $\rho_\ell(\alpha)$ of $F(s,\alpha)$ at $s=s^*_\ell$ is given by
\begin{equation}
\notag
\rho_\ell(\alpha) = \frac{d_\ell(F)}{d} \frac{\omega_F}{\sqrt{2\pi}} e^{-i\frac{\pi}{2}(\xi_F+ds_\ell^*)} \left(\frac{q^{1/d}}{2\pi d}\right)^{\frac{d}{2}-ds_\ell^*} \frac{\overline{a(n_\alpha)}}{n_\alpha^{1-s_\ell^*}}.
\end{equation}

\smallskip
If $F\in\S^\sharp$ has $d=2$, $q=1$ and is normalized, the above quantities simplify as follows:
\begin{equation}
\label{2-12}
\begin{split}
\theta_F=0, \quad \text{Spec}(F) &= \{2\sqrt{m}: m\in\NN \ \text{with} \ a(m)\neq 0\}, \quad n_\alpha=\alpha^2/4, \quad s^*_\ell=s_\ell = \frac34-\frac{\ell}{2}, \\
h_F(s)& \approx \frac{\omega_F}{\sqrt{2\pi}} (4\pi)^{2s-1} \sum_{\ell=0}^\infty d_\ell(F) \Gamma(3/2-2s-\ell), \qquad d_0(F)=1, \\
&\hskip-.3cm \rho_\ell(\alpha) = \frac{e^{i\pi/4}\overline{a(\alpha^2/4)}}{\sqrt{\alpha}}  d_\ell(F) (-2\pi i)^{-\ell} \alpha^{-\ell} \ \ \text{for $\alpha\in$ Spec$(F)$}.
\end{split}
\end{equation}
Note that we used Lemma 4.1 in Section 4 to get the above simplified expression of $\rho_\ell(\alpha)$. Note also that Spec$(F)$ is an infinite set. Indeed, $F$ cannot be a Dirichlet polynomial since, by the functional equation, its Lindel\"of function $\mu_F(\si)$ is positive for negative values of $\si$.

\medskip
\section{Outline of the proof}

\smallskip
In Section 4 we start the proof with a finer investigation of the properties of the structural invariants $d_\ell(F)$ with $\ell\geq0$, introduced in \cite{Ka-Pe/revisited}, in the case of normalized functions $F\in\S^\sharp$ of degree 2 and conductor 1. These invariants do not characterize completely the functional equation of $F$, but contain a good amount of information on it. In particular, the $d_\ell(F)$ essentially determine the $h$-function $h_F(s)$ appearing in the invariant form \eqref{2-10} of the functional equation; see \eqref{2-12}. The main tool in such investigation is a fully explicit version of a special case of the transformation formula for nonlinear twists studied in \cite{Ka-Pe/2011},\cite{Ka-Pe/2015},\cite{Ka-Pe/2016b},\cite{Ka-Pe/2017}; see Lemma 4.2. From this formula we derive an interesting result in itself, namely: for every integer  $N\geq 2$ the invariants $d_\ell(F)$ with $\ell\leq N$ lie on the algebraic variety defined by $Q_N(X_0,\dots,X_N)=0$, where $Q_N$ are certain {\it universal} quadratic forms, i.e. independent of $F$; see Proposition 4.1. Since $d_0(F)=1$, an inductive argument immediately shows that all $d_\ell(F)$ with $\ell\geq 2$ are determined by $d_1(F)$ by means of a procedure not depending on $F$. Then we derive the value of $d_1(F)$ in terms of $H$-invariants, which are easy to compute from the data of the functional equation of $F$. Actually, recalling \eqref{1-1} and \eqref{2-5}, it turns out that 
\begin{equation}
\label{3-1}
d_1(F) = \chi_F - \frac{1}{8};
\end{equation}
see Lemma 4.3. Moreover, $d_1(F)$ is real-valued in the present case.

\smallskip
Next, in Section 5 we introduce the {\it virtual $\ga$-factors}
\begin{equation}
\label{3-2}
\gamma(s) = 
\begin{cases}
(2\pi)^{-s} \Gamma(s+\mu) &\text{with $\mu>0$} \\
\pi^{-s} \Gamma\big(\frac{s+\epsilon+i\kappa}{2}\big) \Gamma\big(\frac{s+\epsilon-i\kappa}{2}\big) &\text{with $\epsilon\in\{0,1\}$ and $\kappa\geq 0$,}
\end{cases}
\end{equation}
respectively of Hecke and Maass type, and consider the analogues $h_\ga(s)$ of the $h$-function, $S_\ga(s)$ of the $S$-function in \eqref{5-14} and $d_\ell(\ga)$ of the structural invariants, see \eqref{5-8}. Although not every virtual $\ga$-factor corresponds to an existing $L$-function, the invariants $d_\ell(\ga)$ satisfy the same properties of the $d_\ell(F)$; see Lemma 5.2. Moreover, since
\begin{equation}
\label{3-3}
\chi_\ga:= \xi_\gamma + H_\gamma(2) + \frac{2}{3}= 
\begin{cases}
2\mu^2 \\
- 2\kappa^2
\end{cases}
\end{equation}
with $\mu$ and $\kappa$ as in \eqref{3-2}, thanks to the analog of \eqref{3-1} the set $\{d_1(\ga):   \text{$\ga$ virtual $\ga$-factor}\}$ coincides with $\RR$. Hence to any $F$ we associate a unique virtual $\ga$-factor such that 
\begin{equation}
\label{3-4}
d_1(F)=d_1(\ga).
\end{equation}
Thus $d_\ell(F)=d_\ell(\ga)$ for every $\ell\geq0$ by the above reported properties of the structural invariants, and hence $h_F(s) = \omega_F h_\ga(s)$; see Proposition 5.1. As a consequence, $F$ satisfies the functional equation 
\begin{equation}
\label{3-5}
\gamma(s) F(s) = \omega_F R(s) \gamma(1-s) F(1-s), \qquad R(s) = \frac{S_F(s)}{S_\ga(s)},
\end{equation}
where $\gamma$ is the virtual $\ga$-factor associated with $F$; see Corollary 5.1. Now we observe that if the function $R(s)$ is constant, then \eqref{3-5} becomes a functional equation of Hecke or Maass type, thus Theorem 1.1 follows at once from classical converse theorems, see Lemma 5.1, thanks to \eqref{3-1},\eqref{3-3} and \eqref{3-4}. Moreover, in Lemma 5.3 we show that if $R(s)$ is not constant, then the function $S_F(s)$, rewritten as
\[
S_F(s) = \sum_{j=0}^N a_j e^{i\pi\omega_js} \quad \text{with $a_j\neq0$ and $\omega_j$ strictly increasing,}
\]
has $N\geq3$ and hence contains the term
\begin{equation}
\label{3-6}
a_{N-1}e^{i\pi\omega_{N-1}s} \quad \text{with $\omega_{N-1}>0$.}
\end{equation}

\smallskip
The last step, i.e. proving that $R(s)$ is constant, starts in Section 6 and involves a non-standard use of certain {\it period functions}, in the sense of Lewis-Zagier \cite{Le-Za/2001}. Indeed, given $F$ and its virtual $\ga$-factor, in Section 6 we consider the Fourier series
\[
f(z) = \sum_{n=1}^\infty a(n) n^\lambda e(nz), \qquad z\in\HH,
\]
where $\la=\mu$ or $\la=i\kappa$ according to \eqref{3-2}, and the period function
\[
\psi(z) = f(z) - z^{-2\la-1} f(-1/z).
\]
Then we proceed to the analysis of the function $f(z)$, involving the properties of certain Mittag-Leffler functions and of the three-term functional equation \eqref{6-23} satisfied by $\psi(z)$. Such an analysis requires different arguments according to the type of the associated virtual $\ga$-factor, the Hecke case being more delicate than the Maass case. For this reason we add the suffixes $H$ and $M$ to distinguish the two cases. The results of the above analysis are contained in Propositions 6.1 and 6.2, where we show that
\[
\psi(z) = Q_H(z) + H_1(z) \qquad \text{and} \qquad \psi(z) = Q_M(z) + H_2(z),
\]
respectively, where all the involved functions are holomorphic for $|\arg(z)|<\pi$ and $Q_H(z)$, $Q_M(z)$ are certain integrals. For example, in the Maass case we have
\[
Q_M(z) = \frac{\omega_F}{2\pi i} \int_{(c_1)} \big(R(s-\la)-1) \Gamma(s) \frac{\ga(1-s+\la)}{\ga(s-\la)} F(1-s+\la) (-2\pi iz)^{-s} \d s
\]
with $0<c_1<1$.

\smallskip
Finally, in Section 7 we analyse the integrals $Q_H(z)$ and $Q_M(z)$. Our aim is to rebuild the function $f(z)$, and to exploit the range where such integrals are holomorphic in order to continue $f(z)$ to an angular region below the real axis. In this process we assume that $R(s)$ is not constant, and hence $N\geq 3$ and \eqref{3-6} holds. This is crucial to show that, roughly, $Q_H(z)$ and $Q_M(z)$ transform to 
\[
f\big(e^{i\pi (1-\omega_{N-1})}  z\big) +H_3(z), 
\]
where $H_3(z)$ is holomorphic in a suitable angular region. After a change of variable we obtain that $f(z)$ is holomorphic for $-\rho\pi<\arg(z)<\pi$ with some $\rho>0$. But $f(z)$ is 1-periodic and hence it is an entire function; this leads immediately to a contradiction. Thus $R(s)$ is constant and Theorem 1.1 follows.

\medskip
\section{Invariants}

\smallskip
We refer to Sections 1 and 2 for definitions and notation used in what follows. We start recalling the basic properties of normalized functions in $\S^\sharp$ with degree $2$ and conductor $1$.

\medskip
{\bf Lemma 4.1.} {\sl Let $F\in\S^\sharp$ be of degree $2$ and conductor $1$. Then $F$ can be normalized. If $F$ is normalized then its Dirichlet coefficients are real, $\xi_F$ is an even integer and $\omega_F = -e^{i\frac{\pi}{2}\xi_F}$. Moreover, all its $H$-invariants are real and its $\gamma$-factor satisfies $\overline{\gamma}(s) = \gamma(s)$.}

\medskip
{\it Proof.} If $F$ has a pole at $s=1$ then $\theta_F=0$, see Lemma 4.1 of \cite{Ka-Pe/2015}, while an entire $F$ can always be shifted by a purely imaginary quantity to get $\theta_F=0$. Clearly, such a shift does not change degree and conductor. Hence, if necessary, every $F\in\S^\sharp$ with $d=2$ and $q=1$ can be normalized by a suitable shift and multiplication by a constant. The second statement follows from equation (5) and the theorem in \cite{Ka-Pe/2017}, observing that $\omega^*_F$ there is denoted here by $\omega_F$ and $H_F(1)=\xi_F$, see \eqref{2-5}. The assertion about $H$-invariants follows since the $\mu$-data of the conjugate function $\overline{F}$ are $\overline{\mu}_j$, $j=1,\dots,r$, thus $H_{\overline{F}}(n)=\overline{H_F(n)}$ for every $n\geq0$ and hence $H_F(n)\in\RR$ in our case. This also explains why $\overline{\gamma}(s) = \gamma(s)$. \qed

\medskip
Since every $F\in\S^\sharp$ with $d=2$ and $q=1$ can be normalized, from now on we always consider normalized functions, as in Theorem 1.1.

\medskip
A major source of information about invariants comes from the transformation formula for nonlinear twists, which we studied in \cite{Ka-Pe/2011},\cite{Ka-Pe/2015},\cite{Ka-Pe/2016b},\cite{Ka-Pe/2017}. Indeed, roughly speaking, the transformation formula gives different outputs if the same nonlinear twist  of a function $F\in\S^\sharp$ is written in formally different ways; this phenomenon imposes several constraints on the invariants. Actually, some of the results in Lemma 4.1 were obtained along these lines, and to proceed further we need a closer analysis of the transformation formula. This is done, in a special case relevant for our purposes, in the next lemma. For $\alpha>0$ we rewrite the standard twist $F(s,\alpha)$ in the form
\[
F(s;f) = \sum_{n=1}^\infty \frac{a(n)}{n^s} e(-f(n,\alpha)), \quad f(n,\alpha) = n + \alpha\sqrt{n}.
\]
For any integer $m\geq0$ we define the polynomials in the $(s,\alpha)$-variables
\begin{equation}
\label{4-1}
\begin{split}
W_m(s,\alpha) &=\mathop{\sum_{\nu=0}^\infty \sum_{\mu=0}^\infty \sum_{k=3\nu}^\infty  \sum_{\ell=0}^{\infty} \sum_{h=0}^\infty}_{\substack{2|(\mu+k) \\ -2\nu+\mu+k+2\ell+h=m}}  A(\nu,\mu,k,\ell,h) \\
&\times {-\frac14-s-\frac{\ell}{2} \choose \mu}  {\frac12-2s +2\nu-\mu-k-\ell \choose h} d_\ell(F) \alpha^h,
 \end{split}
\end{equation}
where $d_\ell(F)$ are the structural invariants of $F$,
\begin{equation}
\label{4-2}
A(\nu,\mu,k,\ell,h) = \frac{1}{\sqrt{\pi} \nu!} \big(-\frac{2i}{\pi}\big)^{\frac{\mu+k}{2}} \big(-\frac12\big)^h (4\pi)^{\nu-\ell} a_{k,\nu} \Gamma\big(\frac{\mu+k+1}{2}\big) i^{\nu+\ell}
\end{equation}
and the coefficients $a_{k,\nu}$ are defined by the expansion
\begin{equation}
\label{4-3}
\Big(\sum_{k=3}^\infty {1/2 \choose k} \xi^k\Big)^\nu = \sum_{k=3\nu}^\infty a_{k,\nu} \xi^k.
\end{equation}
Note that $W_m(s,\alpha)$ is a polynomial since $k\geq 3\nu$ and hence equation $-2\nu+\mu+k+2\ell+h=m$ in \eqref{4-1} has only finitely many solutions for every $m$.

\medskip
{\bf Lemma 4.2.} {\sl Let $F$ be as in Theorem $1.1$ and $\alpha\in$\hskip.1cm{\rm Spec}$(F)$. Then, with the above notation, for every integer $M\geq0$ we have
\[
F(s;f) =  \sum_{m=0}^{M} W_m(s,\alpha) F\big(s+\frac{m}{2},\alpha\big) + H_{M}(s,\alpha),
\]
where $W_0(s,\alpha)\equiv 1$ and the function $ H_{M}(s,\alpha)$ is holomorphic for $\si>-(M-1)/2$.}

\medskip
{\it Proof.} We follow closely the arguments in \cite{Ka-Pe/2011} and \cite{Ka-Pe/2016b}, see also \cite{Ka-Pe/2017}, giving details only at the places which are important for this paper. Moreover, since our goal here is to compute explicitly the quantities appearing in the main term of the transformation formula, we proceed in an essentially formal way, largely disregarding the problems about convergence, regularity and error terms. Indeed, for these issues a detailed treatment already appears in the above mentioned papers. Therefore, in order to simplify the presentation we denote generically by $\E$ any quantity not contributing to our explicit computations, and replace equality by $\sim$ whenever a finite sum plus error term is replaced by the full expansion. We also correct some minor inaccuracies in our previous computations. Such inaccuracies were unimportant for the results in \cite{Ka-Pe/2011} and \cite{Ka-Pe/2016b}, while precise computations of the involved quantities are essential here. Moreover, we simply write $d_\ell$ for the structural invariants $d_\ell(F)$.

\smallskip
(a) {\it Set up.} 
Let $\alpha>0$, not necessarily in Spec$(F)$. Since in our case
\[
f(n,\alpha) = n + \alpha\sqrt{n} = n^{\kappa_0} + \alpha n^{\kappa_1},
\]
say, as in Section 2 of \cite{Ka-Pe/2016b} for $X>1$ and $s\in\CC$ we write
\[
F_X(s;f) = \sum_{n=1}^\infty \frac{a(n)}{n^s} e(-f(n,\alpha)) e^{-(n+\sqrt{n})/X}.
\]
Moreover, as in Section 2.1 of \cite{Ka-Pe/2016b}, let $R\geq1$ be such that $3/2 +2R\not\in\NN$, $-R<\si<-R+\delta$ with a small $\delta>0$, $\rho>2(R+1)/3$, $\emptyset\neq \A\subset\{0,1\}$ with cardinality $|\A|$,
\[
\bfw = w_0+\frac12w_1, \quad \d \bfw = \d w_0 \d w_1, \quad z_0=\frac{1}{X} +2\pi i, \quad z_1=\frac{1}{X} + 2\pi i\alpha,
\]
\[
G(\bfw) = \Gamma(w_0)\Gamma(w_1) z_0^{-w_0} z_1^{-w_1}
\]
and let $\bfw_{|\A}$, $\d\bfw_{|\A}$ and $G(\bfw_{|\A})$ be the restriction of $\bfw$, $\d\bfw$ and $G(\bfw)$ to the indices appearing in $\A$, respectively. For $\si>-R$ and $\Re(w_\nu)=\rho$ for $\nu=0,1$ we have $\Re(s+\bfw)>1$, hence by Mellin's transform we get, as in (2.1) of \cite{Ka-Pe/2016b}, that
\[
F_X(s;f) = \frac{1}{(2\pi i)^2} \int_{(\rho)}\int_{(\rho)} F(s+\bfw) G(\bfw) \d\bfw.
\]
Next we shift the lines of integration to $-\eta$ with some $1/2<\eta<3/4$. Writing
\[
\int_{\LL} \d\bfw = \int_{(-\eta)} \int_{(-\eta)} \d\bfw \quad \text{and analogously for} \quad \int_{\LL_{|\A}} \d\bfw_{|\A},
\]
and
\[
I_X(s,\A) = \frac{1}{(2\pi i)^{|\A|}} \int_{\LL_{|\A}} F(s+\bfw_{|\A}) G(\bfw_{|\A}) \d\bfw_{|\A},
\]
as in Lemma 2.1 of \cite{Ka-Pe/2016b} we have that
\begin{equation}
\label{4-4}
F_X(s;f) =  \sum_{\emptyset\neq\A\subset\{0,1\}} I_X(s,\A) + \E.
\end{equation}

\smallskip
Now, recalling that the coefficients of $F$ are real by Lemma 4.1, since $\si-\eta<0$ we apply  inside each integral $I_X(s,\A)$ the functional equation of $F$ in the invariant form
\begin{equation}
\label{4-5}
F(s+\bfw_{|\A}) = h_F(s+\bfw_{|\A}) S_F(s+\bfw_{|\A}) F(1-s-\bfw_{|\A}),
\end{equation}
see \eqref{2-10}, and expand $F(1-s-\bfw_{|\A})$. This way, writing $\bfx_{|\A}= s+\bfw_{|\A}$ as on p.7656 of \cite{Ka-Pe/2016b}, we obtain that
\begin{equation}
\label{4-6}
I_X(s,\A) = \sum_{n=1}^\infty \frac{a(n)}{n^{1-s}}  \frac{1}{(2\pi i)^{|\A|}} \int_{\LL_{|\A}} h_F(\bfx_{|\A}) S_F(\bfx_{|\A}) G(\bfw_{|\A}) n^{\bfw_{|\A}} \d\bfw_{|\A}.
\end{equation}
Note that the right hand side of \eqref{4-6} is formally slightly different from (2.10) of \cite{Ka-Pe/2016b}, as here we use the invariant form \eqref{4-5} of the functional equation. Precisely, comparing with the notation in (2.10) of \cite{Ka-Pe/2016b} we have
\[
S_F(\bfx_{|\A}) = 2^r S(\bfx_{|\A}) \quad \text{and} \quad h_F(\bfx_{|\A}) = \frac{\omega}{(2\pi)^r} Q^{1-2\bfx_{|\A}} \widetilde{G}(\bfx_{|\A}).
\]
However, this slight difference does not affect the subsequent arguments in \cite{Ka-Pe/2011},\cite{Ka-Pe/2016b},\cite{Ka-Pe/2017}. Actually, our present approach simplifies the treatment since we can appeal directly to the asymptotic expansion in \eqref{2-12} of the function $h_F(s)$, involving the structural invariants $d_\ell$. The presence of the invariants $d_\ell$ provides also a conceptual advantage, which turns out to be important later on in the paper.

\smallskip
(b) {\it Expansions.}
From (2.12) of \cite{Ka-Pe/2016b} we have that 
\[
S_F(\bfx_{|\A}) = a_{-N} e^{-i\pi \bfx_{|\A}} + a_{N} e^{i\pi \bfx_{|\A}} + \E.
\]
From the treatment on pages 7657--7659 of \cite{Ka-Pe/2016b} of the $J_X$-integrals, stemming from the term $a_{N} e^{\pi i \bfx_{|\A}}$, we have that in the end their contribution is absorbed by the function $H_M(s,\alpha)$ in the statement of the lemma; see also Lemma 2.2 of \cite{Ka-Pe/2016b}. Thus the term $a_{N} e^{\pi i \bfx_{|\A}}$ does not contribute to the main terms, 
hence we rewrite the above expression for $S_F(\bfx_{|\A})$ as
\[
S_F(\bfx_{|\A}) = a_{-N} e^{-i\pi \bfx_{|\A}} + \E.
\]
Next we recall that $a_{-N} = e^{-i\frac{\pi}{2} \xi_F}$, see Section 3.7 of \cite{Ka-Pe/revisited}, hence $\omega_F a_{-N}=-1$ by Lemma 4.1. Therefore, in view of the expansion of $h_F(s)$ in  \eqref{2-12} we have that equation (2.13) of \cite{Ka-Pe/2016b} is now replaced by
\begin{equation}
\label{4-7}
h_F(\bfx_{|\A}) S_F(\bfx_{|\A}) \sim \frac{-(4\pi)^{2\bfx_{|\A}}}{2(2\pi)^{3/2}} e^{-i\pi\bfx_{|\A}} \sum_{\ell=0}^\infty d_\ell \Gamma(3/2-2\bfx_{|\A}-\ell) + \E.
\end{equation}
Thus, by \eqref{4-6} and \eqref{4-7}, equation (2.15) of \cite{Ka-Pe/2016b} finally becomes
\begin{equation}
\label{4-8}
I_X(s,\A) \sim \frac{-(4\pi)^{2s}}{2(2\pi)^{3/2}} \sum_{\ell=0}^\infty d_\ell  \sum_{n=1}^\infty \frac{a(n)}{n^{1-s}} I_X(s,\A,n,\ell) + \E
\end{equation}
with
\[
I_X(s,\A,n,\ell) =  \frac{1}{(2\pi i)^{|\A|}} \int_{\LL_{|\A}} \Gamma(3/2-2\bfx_{|\A} - \ell\big) e^{-i\pi\bfx_{|\A}} G(\bfw_{|\A}) \big((4\pi)^2n\big)^{\bfw_{|\A}} \d\bfw_{|\A}.
\]
Note that $I_X(s,\A,n,\ell)$ coincides with the definition of $I_X(s,\A,n,\ell)$ in (2.16) of \cite{Ka-Pe/2016b} with $d=2$, $q=1/(4\pi)^2$ and $\theta_F=0$, but not with the definition of $I_X(s,\A,n,\ell)$ on p.1409 of \cite{Ka-Pe/2011}. Indeed, recalling the data of the functional equation of $F$ in Section 2, the latter coincides with
\begin{equation}
\label{4-9}
\big(\frac{4}{\beta}\big)^s I_X(s,\A,n,\ell), \quad \text{where} \quad \beta=\prod_{j=1}^r\lambda_j^{2\lambda_j},
\end{equation}
since $Q^2\beta/d^d=(4\pi)^{-2}$ in the present case, where degree is 2 and conductor is 1.

\smallskip
Now we use the Mellin transform argument leading to equation (2.18) of \cite{Ka-Pe/2016b}, but in order to compute the coefficients $f_\ell$ in such equation we have to follow the arguments on p.1410 of \cite{Ka-Pe/2011}. Recalling the difference in notation pointed out in \eqref{4-9}, by (2.13) of \cite{Ka-Pe/2011} we have that the inverse Mellin transform $\widetilde{I}_X(y)$ of $\big(\frac{4}{\beta}\big)^s I_X(s,\A,n,\ell)$ is in our case
\[
\widetilde{I}_X(y) = \frac12 e^{2i/\sqrt{y\beta}} \Big(\frac{2e^{-i\pi/2}}{\sqrt{y\beta}}\Big)^{\frac{3}{2}-\ell} \prod_{\nu\in\A} \Big(e^{-\frac{z_\nu Q^{2\kappa_\nu}}{n^{\kappa_\nu}y^{\kappa_\nu}}} - 1\Big).
\]
Hence, after a computation needed since the constants $c_\ell$ in (2.13) and (2.14) of \cite{Ka-Pe/2011} are not the same, equation (2.14) of \cite{Ka-Pe/2011} in the present case becomes
\[
\big(\frac{4}{\beta}\big)^s I_X(s,\A,n,\ell) = \big(\frac{4}{\beta}\big)^s c_\ell \int_0^\infty e^{i\sqrt{x}} \prod_{\nu\in\A} \Big(e^{-z_\nu (\frac{x}{(4\pi)^2 n})^{\kappa_\nu}} - 1\Big) x^{-s-\frac{\ell}{2} - \frac{1}{4}} \d x,
\]
where
\begin{equation}
\label{4-10}
c_\ell = \frac12 e^{-\frac{3}{4}\pi i} i^\ell.
\end{equation}
As a consequence, using identity (2.16) of \cite{Ka-Pe/2011}, equation (2.19) of \cite{Ka-Pe/2016b} now takes the form
\begin{equation}
\label{4-11}
\begin{split}
I_X(s,n,\ell) &:= \sum_{\emptyset\neq\A\subset\{0,1\}} I_X(s,\A,n,\ell) \\
&=  c_\ell \int_0^\infty e^{i\sqrt{x}} \big(e^{-\Psi_X(x,n)} e\big(-f\big(\frac{x}{(4\pi)^2 n},\alpha\big)\big)-1\Big) x^{-s-\frac{\ell}{2} - \frac{1}{4}} \d x
\end{split}
\end{equation}
with
\[
\Psi_X(x,n) = \frac{1}{X} \Big(\frac{x}{(4\pi)^2 n} + \sqrt{\frac{x}{(4\pi)^2 n}}\Big).
\]

\smallskip
(c) {\it Saddle point and limit as $X\to\infty$.}
The next step is to apply the saddle point technique developed in Section 2.3 of \cite{Ka-Pe/2011} and \cite{Ka-Pe/2016b}  to the right hand side of \eqref{4-11}. Recalling that $\theta_F=0$ and that the integral on the right hand side of \eqref{4-11} is multiplied by $c_\ell$, from Lemma 2.3 of \cite{Ka-Pe/2016b} we deduce that for $n\geq n_0$ and $n_0$ sufficiently large
\begin{equation}
\label{4-12}
I_X(s,n,\ell) = c_\ell K_X\big(s+ \frac{\ell}{2},n\big) + \E,
\end{equation}
where $K_X(s,\xi)$ is as in (2.21) of \cite{Ka-Pe/2016b}. Hence from \eqref{4-4},\eqref{4-8},\eqref{4-10} and \eqref{4-12} we have that (2.22) of \cite{Ka-Pe/2016b} becomes
\[
F_X(s;f) \sim \frac{e^{i\pi/4}}{4(2\pi)^{3/2}} (4\pi)^{2s} \sum_{\ell=0}^\infty i^\ell d_\ell  \sum_{n\geq n_0} \frac{a(n)}{n^{1-s}} K_X(s+\frac{\ell}{2},n) + \E
\]
Then we pass to the limit as $X\to\infty$ as in Section 2.4 of \cite{Ka-Pe/2016b}, thus getting that
\begin{equation}
\label{4-13}
F(s;f) \sim \frac{e^{i\pi/4}}{4(2\pi)^{3/2}} (4\pi)^{2s} \sum_{\ell=0}^\infty i^\ell d_\ell  \sum_{n\geq n_0} \frac{a(n)}{n^{1-s}} K(s+\frac{\ell}{2},n) + \E,
\end{equation}
where
\[
K(s,n) = \gamma x_0^{\frac{3}{4}-s} \int_{-r}^r e^{i\Phi(z,n)} (1+\gamma\lambda)^{-s-\frac{1}{4}} \d \lambda
\]
and $x_0=x_0(n)$ is the critical point of $\Phi(z,n)$, see Lemma 2.3 of \cite{Ka-Pe/2011},
\[
\Phi(z,n) = z^{\frac12} - 2\pi f\big(\frac{z}{(4\pi)^2n},\alpha\big), \quad z=x_0(1+\gamma\lambda),
\] 
\[
\gamma = 1-i, \quad r= \frac{\log n}{\sqrt{R}}, \quad R=x_0^2\Phi''(x_0,n).
\]
We finally note that the critical point $x_0$ is real and
\begin{equation}
\label{4-14}
\frac{\log n}{\sqrt{n}} \ll r \ll \frac{\log n}{\sqrt{n}},
\end{equation}
see Lemma 2.3 and (2.30) of \cite{Ka-Pe/2011}, respectively. Note also that, with respect to \cite{Ka-Pe/2011} and \cite{Ka-Pe/2016b}, we slightly simplified the notation of the various terms in $K(s,\xi)$.

\smallskip
(d) {\it Computation of $K(s,n)$.}
Now we proceed to the novel part of the proof of the lemma, namely a detailed computation of the function $K(s,n)$. From now on we assume that $\alpha\in$ Spec$(F)$, thus $\alpha^2/4\in\NN$. Clearly, the critical point $x_0$, solution of the equation $\frac{\partial}{\partial z} \Phi(z,n) =0$, in our case satisfies
\begin{equation}
\label{4-15}
\sqrt{x_0} = 4\pi n \Delta_n \quad \text{with} \quad \Delta_n = 1-\frac{\alpha}{2\sqrt{n}},
\end{equation}
hence
\[
\Phi(x_0,n) = 2\pi\big(n-\alpha\sqrt{n}+\frac{\alpha^2}{4}\big) \quad \text{and} \quad e^{i\Phi(x_0,n)} = e(-\alpha\sqrt{n}).
\]
Therefore
\begin{equation}
\label{4-16}
\begin{split}
K(s,n) &= \gamma x_0^{\frac{3}{4}-s} e^{i\Phi(x_0,n)} \int_{-r}^r e^{i(\Phi(z,n)-\Phi(x_0,n))} (1+\gamma\lambda)^{-s-\frac{1}{4}} \d \lambda \\
&= \gamma e(-\alpha\sqrt{n}) (4\pi n\Delta_n)^{\frac{3}{2}-2s} I(s,n),
\end{split}
\end{equation}
where
\[
I(s,n) = \int_{-r}^r e^{i(\Phi(z,n)-\Phi(x_0,n))} (1+\gamma\lambda)^{-s-\frac{1}{4}} \d \lambda.
\]
But
\[
\begin{split}
\Phi(z,n)-\Phi(x_0,n) &= \sqrt{z} \Delta_n - \frac{z}{8\pi n} - \sqrt{x_0} \Delta_n + \frac{x_0}{8\pi n} \\
&= 4\pi n\Delta_n^2 \big( (1+\gamma\lambda)^{1/2} -1 - \frac{\gamma\lambda}{2} \big),
\end{split}
\]
hence writing
\begin{equation}
\label{4-17}
\phi(x) = (1+x)^{1/2} -1 - \half x + {\textstyle{1\over 8}} x^2 = \sum_{k=3}^\infty {1/2 \choose k} x^k
\end{equation}
and recalling that $\gamma^2=-2i$ we have that
\begin{equation}
\label{4-18}
\begin{split}
I(s,n) &= \int_{-r}^r e^{4\pi in\Delta_n^2(\phi(\gamma\lambda) + \frac{1}{4}i\lambda^2)} (1+\gamma\lambda)^{-s-\frac{1}{4}} \d \lambda \\
&=  \int_{-r}^r e^{-\pi n\Delta_n^2\lambda^2} e\big(2n\Delta_n^2 \phi(\gamma\lambda)\big) (1+\gamma\lambda)^{-s-\frac{1}{4}} \d \lambda.
\end{split}
\end{equation}

\smallskip
Since $r=o(1)$ as $n\to\infty$ by \eqref{4-14}, for $n\geq n_0$ we expand the complex exponential and the function $(1+\gamma\lambda)^{-s-\frac{1}{4}}$ in the above integral, and then we replace $\phi(\gamma\lambda)$ by its power series. Hence we have
\[
e\big(2n\Delta_n^2 \phi(\gamma\lambda)\big) (1+\gamma\lambda)^{-s-\frac{1}{4}} = \sum_{\nu=0}^\infty \frac{\big(4\pi in\Delta_n^2\phi(\gamma\lambda)\big)^\nu}{\nu!} \sum_{\mu=0}^\infty {-s-\frac14 \choose \mu} (\gamma\lambda)^\mu
\]
and by \eqref{4-17} and \eqref{4-3} we write
\[
\phi(\gamma\lambda)^\nu = \sum_{k=3\nu}^\infty a_{k,\nu} (\gamma\lambda)^k \quad \text{with} \quad a_{0,0}=1.
\]
Thus, thanks to the symmetry of the integral, \eqref{4-18} becomes
\begin{equation}
\label{4-19}
I(s,n) =  2\mathop{\sum_{\nu=0}^\infty \sum_{\mu=0}^\infty  \sum_{k=3\nu}^\infty}_{2|(\mu+k)} \frac{(4\pi in\Delta_n^2)^\nu}{\nu!} {-s-\frac14 \choose \mu} a_{k,\nu} \gamma^{\mu+k} \int_0^r e^{-\pi n\Delta_n^2\lambda^2} \lambda^{\mu+k} \d\lambda.
\end{equation}
But, first by the change of variable $\sqrt{n}\Delta_n\lambda=u$, then completing the integral to $(0,\infty)$ and finally by the change of variable $\pi u^2=\xi$, in view of \eqref{4-14} we have
\[
\begin{split}
\int_0^r e^{-\pi n\Delta_n^2\lambda^2} \lambda^{\mu+k} \d\lambda &= \frac{1}{(\sqrt{n}\Delta_n)^{\mu+k+1}} \int_0^\infty e^{-\pi u^2} u^{\mu+k} \d u + \E  \\
&= \frac{1}{2(\sqrt{\pi n}\Delta_n)^{\mu+k+1}} \int_0^\infty e^{-\xi} \xi^{\frac{\mu+k-1}{2}} \d\xi + \E \\
& =  \frac{\Gamma\big(\frac{\mu+k+1}{2}\big)}{2(\sqrt{\pi n}\Delta_n)^{\mu+k+1}} +\E.
\end{split}
\]
Hence \eqref{4-19} becomes
\begin{equation}
\label{4-20}
I(s,n) =   \mathop{\sum_{\nu=0}^\infty \sum_{\mu=0}^\infty  \sum_{k=3\nu}^\infty}_{2|(\mu+k)} \frac{(4\pi i)^\nu}{\nu!} a_{k,\nu} \gamma^{\mu+k} \frac{\Gamma\big(\frac{\mu+k+1}{2}\big)}{\pi^{\frac{\mu+k+1}{2}}} {-s-\frac14 \choose \mu} \frac{\Delta_n^{2\nu-\mu-k-1}}{n^{-\frac{{2\nu-\mu-k-1}}{2}}} + \E,
\end{equation}
and from \eqref{4-16} and \eqref{4-20} we obtain that
\begin{equation}
\label{4-21}
\begin{split}
K(s,n) &= \gamma e(-\alpha\sqrt{n}) (4\pi)^{\frac{3}{2}-2s} (n\Delta_n)^{\frac{3}{2}-2s} \\
&\hskip-1cm \times  \mathop{\sum_{\nu=0}^\infty \sum_{\mu=0}^\infty  \sum_{k=3\nu}^\infty}_{2|(\mu+k)} \frac{(4\pi i)^\nu}{\nu!} a_{k,\nu} \gamma^{\mu+k} \frac{\Gamma\big(\frac{\mu+k+1}{2}\big)}{\pi^{\frac{\mu+k+1}{2}}} {-s-\frac14 \choose \mu} \frac{\Delta_n^{2\nu-\mu-k-1}}{n^{-\frac{{2\nu-\mu-k-1}}{2}}} + \E \\
&= \mathop{\sum_{\nu=0}^\infty \sum_{\mu=0}^\infty  \sum_{k=3\nu}^\infty}_{2|(\mu+k)} b_{\nu,\mu,k}(s)  \frac{\Delta_n^{2\nu-\mu-k+\frac12-2s}}{n^{2s-1-\frac{{2\nu-\mu-k}}{2}}} e(-\alpha\sqrt{n})+ \E,
\end{split}
\end{equation}
where
\[
b_{\nu,\mu,k}(s) = (4\pi)^{\frac{3}{2}-2s} \frac{(4\pi i)^\nu}{\nu!} a_{k,\nu} \gamma^{\mu+k+1} \frac{\Gamma\big(\frac{\mu+k+1}{2}\big)}{\pi^{\frac{\mu+k+1}{2}}}  {-s-\frac14 \choose \mu}.
\]
Note that $b_{\nu,\mu,k}(s)$ is independent of $\alpha$, $n$ and the data of the functional equation of $F$.

\smallskip
Recalling the definition of $\Delta_n$ in \eqref{4-15} we have the following expansion
\[
\Delta_n^{2\nu-\mu-k+\frac12-2s} = \sum_{h=0}^\infty {2\nu-\mu-k+\frac12-2s \choose h} \frac{(-1)^h \alpha^h}{2^h n^{h/2}},
\]
hence from \eqref{4-13} and \eqref{4-21} (with $s$ replaced by $s+\ell/2$) we have that
\begin{equation}
\label{4-22}
F(s;f) \sim  \mathop{\sum_{\nu=0}^\infty \sum_{\mu=0}^\infty  \sum_{k=3\nu}^\infty \sum_{\ell=0}^\infty \sum_{h=0}^\infty}_{2|(\mu+k)}  c_{\nu,\mu,k,\ell,h}(s) d_\ell  \alpha^h  \sum_{n\geq n_0} \frac{a(n)}{n^{s+\ell -\frac{2\nu-\mu-k-h}{2}}} e(-\alpha\sqrt{n}) + \E
\end{equation}
with
\begin{equation}
\label{4-23}
c_{\nu,\mu,k,\ell,h}(s) = \frac{e^{i\pi/4}}{4(2\pi)^{3/2}} (4\pi)^{2s} i^\ell {2\nu-\mu-k+\frac12-2s-\ell \choose h} \frac{(-1)^h}{2^h} b_{\nu,\mu,k}\big(s+\frac{\ell}{2}\big).
\end{equation}
Note that, again, the coefficients $c_{\nu,\mu,k,\ell,h}(s)$ are independent of $\alpha$, $n$ and the data of the functional equation of $F$. Hence, recalling that $\gamma=1-i$ and rearranging terms, a comparison of \eqref{4-22} and \eqref{4-23} with \eqref{4-1} and \eqref{4-2} shows that \eqref{4-22} can be rewritten as
\[
F(s;f) \sim \sum_{m=0}^\infty W_m(s,\alpha) \sum_{n\geq n_0} \frac{a(n)}{n^{s+\frac{m}{2}}} e(-\alpha\sqrt{n}) + \E,
\]
where $W_m(s,\alpha)$ is defined by \eqref{4-1}. Thus, adding and subtracting the terms with $1\leq n<n_0$, recalling the definition of the standard twist of $F$ we finally get
\begin{equation}
\label{4-24}
F(s;f) \sim  \sum_{m=0}^\infty W_m(s,\alpha) F\big(s+\frac{m}{2},\alpha\big) + \E.
\end{equation}
Moreover, from \eqref{4-1} we have that
\begin{equation}
\label{4-25}
W_0(s,\alpha) = A(0,0,0,0,0) =  \frac{\Gamma(1/2)}{\sqrt{\pi}} = 1.
\end{equation}

\smallskip
Finally, as in the general case treated in \cite{Ka-Pe/2011} and \cite{Ka-Pe/2016b}, the explicit meaning of \eqref{4-24} is that for every $M\geq0$ there exists $K=K(M)\to\infty$ as $M\to\infty$ such that
\begin{equation}
\label{4-26}
F(s;f) =  \sum_{m=0}^{M} W_m(s,\alpha) F\big(s+\frac{m}{2},\alpha\big) + H_{M}(s,\alpha)
\end{equation}
holds with a function $H_{M}(s,\alpha)$ holomorphic for $\si>-K$. Given $M$, if $K(M)<(M-1)/2$ we choose $M^*>M$ such that $K(M^*)\geq (M-1)/2$, thus getting that \eqref{4-26} holds with $M^*$ in place of $M$ and $-(M-1)/2$ in place of $-K$. But for $M<m\leq M^*$ and $\si>-(M-1)/2$ we have that $\Re(s+m/2) > (1-M)/2 + (M+1)/2 = 1$, hence
\[
\text{$W_m(s,\alpha) F(s+\frac{m}{2},\alpha)$ is holomorphic for $\si>(1-M)/2$.}
\]
The lemma now follows in view of \eqref{4-25} and \eqref{4-26}, bringing the holomorphic terms with $M<m\leq M^*$ to the right hand side. \qed

\medskip
The previous lemma is the first step in the proof of the main result of this section.

\medskip
{\bf Proposition 4.1.} {\sl For every $N\geq2$ there exists a quadratic form 
\[
Q_N(X_0,\dots,X_N) = \sum_{\substack{\ell,h\geq 0 \\ \ell+h\leq N}} \alpha_{\ell,h} X_\ell X_h \quad \text{with $\alpha_{\ell,h}\in\RR$ and $\alpha_{0,N}+a_{N,0}=1$},
\]
such that for every $F$ as in Theorem $1.1$ we have
\[
Q_N\big(d_0(F),\dots,d_N(F)\big) = 0.
\]
Before normalization of coefficients, such a quadratic form is given by \eqref{4-34} and \eqref{4-29} below.}

\medskip
{\it Proof.} Again we write $d_\ell$ for $d_\ell(F)$. Thanks to the periodicity of the complex exponential we have that the nonlinear twist $F(s;f)$ defined before Lemma 4.2 coincides with the standard twist $F(s,\alpha)$, hence Lemma 4.2 implies that for $\alpha\in$ Spec$(F)$ and $M\geq1$
\[
\sum_{m=1}^M W_m(s,\alpha) F\big(s+\frac{m}{2},\alpha\big) \quad \text{is holomorphic for $\si>(1-M)/2$.}
\]
Therefore, the residue of this function at $s=s_M=3/4 - M/2$ is zero. Thus, recalling the polar structure of the standard twist, see Section 2.2, we immediately deduce that
\[
\sum_{m=1}^M W_m(s_M,\alpha) \rho_{M-m}(\alpha)=0.
\]
Hence, thanks to the expression of $\rho_\ell(\alpha)$ in \eqref{2-12} and recalling that $a(\alpha^2/4)\neq0$, we finally obtain that for every $\alpha\in$ Spec$(F)$ and $M\geq1$
\begin{equation}
\label{4-27}
\sum_{m=1}^M W_m(s_M,\alpha) (-2\pi i)^m d_{M-m} \alpha^m=0.
\end{equation}

\smallskip
We set $k=3\nu+b$ with $b\geq0$ in \eqref{4-1} and observe that the range of variation of $2\ell+h$ is
\[
0\leq 2\ell + h \leq m-(\nu+\mu+b) \leq m.
\]
Then, rearranging terms, we rewrite $W_m(s,\alpha)$ as
\begin{equation}
\label{4-28}
\begin{split}
W_m(s,\alpha) &=\mathop{\sum_{\nu=0}^\infty \sum_{\mu=0}^\infty \sum_{b=0}^\infty  \sum_{\ell=0}^{\infty} \sum_{h=0}^\infty}_{\substack{2|(\nu+\mu+b) \\ \nu+\mu+b+2\ell+h=m}}  A(\nu,\mu,3\nu+b,\ell,h) \\
&\hskip1cm \times {-\frac14-s-\frac{\ell}{2} \choose \mu}  {\frac12-2s -\nu-\mu-b-\ell \choose h} d_\ell(F) \alpha^h \\
&= \mathop{\sum_{\ell=0}^\infty \sum_{h=0}^\infty}_{2\ell+h \leq m} B_m(s,\ell,h)  d_\ell(F) \alpha^h,
\end{split}
\end{equation}
where
\begin{equation}
\label{4-29}
\begin{split}
B_m(s,\ell,h)  &=\mathop{\sum_{\nu=0}^\infty \sum_{\mu=0}^\infty \sum_{b=0}^\infty}_{\substack{2|(\nu+\mu+b) \\ \nu+\mu+b=m-(2\ell+h)}}  A(\nu,\mu,3\nu+b,\ell,h) \\
& \times {-\frac14-s-\frac{\ell}{2} \choose \mu}  {\frac12-2s -\nu-\mu-b-\ell \choose h}.
\end{split}
\end{equation}
Thus from \eqref{4-27} and \eqref{4-28} we have, still for every $\alpha\in$ Spec$(F)$ and $M\geq1$, that
\begin{equation}
\label{4-30}
\sum_{m=1}^M (-2\pi i)^m \mathop{\sum_{\ell=0}^\infty \sum_{h=0}^\infty}_{2\ell+h \leq m} B_m(s_M,\ell,h) d_\ell d_{M-m} \alpha^{m+h}=0.
\end{equation}
Next we set $p=m+h$, hence \eqref{4-30} becomes
\begin{equation}
\label{4-31}
\begin{split}
0 &=\sum_{m=1}^M (-2\pi i)^m \mathop{\sum_{\ell=0}^\infty \sum_{p=m}^\infty}_{2\ell+h \leq m} B_m(s_M,\ell,p-m) d_\ell d_{M-m} \alpha^p \\
&= \sum_{1\leq m\leq M} (-2\pi i)^m \sum_{m\leq p \leq 2m} \, \sum_{0\leq \ell \leq m-p/2} B_m(s_M,\ell,p-m) d_\ell d_{M-m} \alpha^p \\
&= \sum_{1 \leq p \leq 2M} \Big(\sum_{p/2 \leq m \leq \min(p,M)}. \sum_{0\leq \ell \leq m-p/2}  (-2\pi i)^m B_m(s_M,\ell,p-m) d_\ell d_{M-m} \Big) \alpha^p
\end{split}
\end{equation}
for every $\alpha\in$ Spec$(F)$ and $M\geq1$. 

\smallskip
But we know that Spec$(F)$ is an infinite set, see at the end of Section 2, thus \eqref{4-31} implies that for every $M\geq1$ and $1\leq p\leq 2M$
\begin{equation}
\label{4-32}
\sum_{p/2 \leq m \leq \min(p,M)}. \sum_{0\leq \ell \leq m-p/2}  (-2\pi i)^m B_m(s_M,\ell,p-m) d_\ell d_{M-m} =0.
\end{equation}
Setting further $h=M-m$, and hence $m=M-h$, we see that condition $p/2 \leq m \leq \min(p,M)$ in \eqref{4-32}  is equivalent to $M-\min(p,M) \leq h \leq M-p/2$ and hence to $\max(M-p,0)\leq h \leq M-p/2$. Therefore we rewrite \eqref{4-32} as
\begin{equation}
\label{4-33}
\sum_{\max(M-p,0) \leq h\leq M-p/2} \ \sum_{0\leq \ell \leq M-h-p/2}  (-2\pi i)^{-h} B_{M-h}(s_M,\ell,p+h-M) d_\ell d_h =0
\end{equation}
for every $M\geq1$ and $1\leq p\leq 2M$. Finally, we choose $p=M$ and $M=2N$ in \eqref{4-33}, thus getting (after reversing summation) that for $N\geq1$
\begin{equation}
\label{4-34}
\widetilde{Q}_N(d_0,\dots,d_N) = \sum_{0 \leq \ell\leq N} \ \sum_{0\leq h \leq N-\ell}  (-2\pi i)^{-h} B_{2N-h}(s_{2N},\ell,h) d_\ell d_h =0.
\end{equation}
This is, before normalization, the quadratic form announced in the proposition. 

\smallskip
Hence it remains to check our assertions about the quadratic form $\widetilde{Q}_N(X_0,\dots,X_N)$, namely that its coefficients are independent of $F$ and real  after normalization, and for $N\geq2$
\[
(-2\pi i)^{-N} B_N(s_{2N},0,N) + B_{2N}(s_{2N},N,0) \neq 0.
\]
The first assertion is immediate since the coefficients $A(\nu,\mu,k,\ell,h)$ in \eqref{4-2}, and hence the polynomials $B_m(s,\ell,h)$ in \eqref{4-29}, are independent of $F$, and so is also $s_{2N}$ by \eqref{2-6}. Moreover, in view of \eqref{4-2} and \eqref{4-29}, the coefficients $(-2\pi i)^{-h} B_{2N-h}(s_{2N},\ell,h)$ are sums of type
\[
\mathop{\sum_{\nu=0}^\infty \sum_{\mu=0}^\infty \sum_{b=0}^\infty}_{\substack{ \nu+\mu+b=2(N-\ell-h)}}  i^{h +\frac{3\nu +\mu+b}{2} +\nu +\ell} c_N(\nu,\mu,b,\ell,h) = i^N \mathop{\sum_{\nu=0}^\infty \sum_{\mu=0}^\infty \sum_{b=0}^\infty}_{\substack{ \nu+\mu+b=2(N-\ell-h)}}  (-1)^\nu c_N(\nu,\mu,b,\ell,h)
\]
with $c_N(\nu,\mu,b,\ell,h) \in\RR$, so the second assertion follows. Finally, again  from \eqref{4-2} and \eqref{4-29} we have that
\[
B_N(s_{2N},0,N) = A(0,0,0,0,N) {\frac12 -2s_{2N} \choose N} = (-2)^{-N} {2N-1\choose N}
\]
and
\[
B_{2N}(s_{2N},N,0) = A(0,0,0,N,0) = \frac{1}{\sqrt{\pi}} (4\pi)^{-N} \Gamma(1/2) i^N = (-4\pi i)^{-N},
\]
thus
\[
(-2\pi i)^{-N} B_N(s_{2N},0,N) + B_{2N}(s_{2N},N,0) = (-4\pi i)^{-N} \Big(1+ (-1)^N{2N-1\choose N}\Big) \neq0
\]
for $N\geq2$, as required. The proposition now follows dividing $\widetilde{Q}_N$ by the right hand side of the last equation. \qed

\medskip
As an immediate consequence of Proposition 4.1 we explicitly record the following corollary, obtained by a trivial induction since $d_0(F)=1$, see \eqref{2-12}.

\medskip
{\bf Corollary 4.1.} {\sl Let $F$ be as in Theorem $1.1$. Then the value of any $d_\ell(F)$ with $\ell\geq2$ is determined by the value of $d_1(F)$ by a recursive algorithm independent of $F$.}

\medskip
Finally we express the structural invariant $d_1(F)$ in terms of the $H$-invariants, which are easy to compute by means of definition \eqref{2-4}.

\medskip
{\bf Lemma 4.3.} {\sl Let $F$ be as in Theorem $1.1$. Then, recalling \eqref{1-1}, we have that}
\[
d_1(F) =  \chi_F - \frac{1}{8}.
\]

\medskip
{\it Proof.} We briefly sketch the formal computations leading to the result; actually, these computations lead to an expression of any $d_\ell(F)$ in terms of $H$-invariants, see \eqref{4-39}, which can be used to calculate explicitly any given $d_\ell(F)$.

\smallskip
Using \eqref{2-3} and Stirling's formula, see Corollary 6.2 of \cite{Vio/2016}, adopting the notation $\sim$ in the proof of Lemma 4.2 we get
\[
\begin{split}
\log h_F(s) &- \log\Gamma(3/2-2s) \sim \log\omega_F - \log\sqrt{2\pi} + (2s-1)\log(4\pi) \\
&+ \sum_{\nu=1}^\infty \frac{(-1)^{\nu+1}}{\nu(\nu+1)} \Big\{\sum_{j=1}^r\Big(\frac{B_{\nu+1}(\lambda_j+\overline{\mu}_j)}{(-\lambda_j s)^\nu} + \frac{B_{\nu+1}(1-\mu_j)}{(-\lambda_j s)^\nu}\Big) - \frac{B_{\nu+1}(3/2)}{(-2s)^\nu}\Big\} \\
&\hskip2.5cm \sim \log\Big(\frac{\omega_F}{\sqrt{2\pi}} \big(\frac{1}{4\pi})^{1-2s}\Big) + \sum_{\nu=1}^\infty \frac{r_\nu(F)}{\nu(\nu+1)}\frac{1}{s^\nu},
\end{split}
\]
say, since $d=2$, $q=1$ and $\theta_F=0$. Moreover, thanks to \eqref{2-4} and to the formulae
\[
B_{\nu+1}(\lambda_j+\overline{\mu}_j) = \sum_{k=0}^{\nu+1} {\nu+1\choose k} B_k(\overline{\mu}_j) \lambda_j^{\nu+1-k} \quad \text{and} \quad B_{\nu+1}(1-\mu_j) = (-1)^{\nu+1} B_{\nu+1}(\mu_j),
\]
see (4.23) and (4.27) of \cite{Vio/2016}, we have that
\begin{equation}
\label{4-35}
r_\nu(F) = \frac{B_{\nu+1}(3/2)}{2^\nu} - \frac12 \Big\{\sum_{k=0}^{\nu+1} {\nu+1\choose k} \overline{H_F(k)} + (-1)^{\nu+1} H_F(\nu+1)\Big\}.
\end{equation}
Thus
\begin{equation}
\label{4-36}
h_F(s) \sim \frac{\omega_F}{\sqrt{2\pi}} \big(\frac{1}{4\pi})^{1-2s} \Gamma(3/2-2s) \exp\Big(\sum_{\nu=1}^\infty \frac{r_\nu(F)}{\nu(\nu+1)}\frac{1}{s^\nu}\Big).
\end{equation}

\smallskip
Next, expanding the exponential we obtain
\[
\exp\Big(\sum_{\nu=1}^\infty \frac{r_\nu(F)}{\nu(\nu+1)}\frac{1}{s^\nu}\Big) \sim 1 + \sum_{m=1}^\infty \frac{1}{m!} \Big(\sum_{\nu=1}^\infty \frac{r_\nu(F)}{\nu(\nu+1)}\frac{1}{s^\nu}\Big)^m \sim 1 + \sum_{h=1}^\infty \frac{V_h(F)}{s^h}
\]
with
\begin{equation}
\label{4-37}
V_h(F) = \sum_{m=1}^h \frac{1}{m!} \mathop{\sum_{\nu_1\geq 1} \dots \sum_{\nu_m\geq1}}_{\substack{\nu_1+\cdots+\nu_m=h}} \prod_{j=1}^m \frac{r_{\nu_j}(F)}{\nu_j(\nu_j+1)} \qquad \text{for} \ h\geq1.
\end{equation}
Moreover, letting $z= 1/2-2s$, for $h\geq 1$ we write
\begin{equation}
\label{4-38}
\frac{1}{s^h} = \sum_{\ell\geq h} \frac{A_{h,\ell}}{z(z-1)\cdots(z-\ell+1)}
\end{equation}
with certain coefficients $A_{h,\ell}\in\RR$. But, by the factorial formula for the $\Gamma$ function, for $\ell\geq1$
\[
\frac{\Gamma(3/2-2s)}{z(z-1)\cdots(z-\ell+1)} = \Gamma(z-\ell+1)=\Gamma(3/2-2s-\ell),
\]
hence substituting into \eqref{4-36} and comparing with \eqref{2-12} we finally obtain that for $\ell\geq1$
\begin{equation}
\label{4-39}
d_\ell(F) =  \sum_{h=1}^\ell A_{h,\ell} V_h(F).
\end{equation}
In particular, \eqref{4-37} and \eqref{4-38} imply that for $\ell=1$ we have $A_{1,1}=-2$ and $V_1(F)=r_1(F)/2$, hence by \eqref{4-39}, \eqref{4-35} and Lemma 4.1
\[
d_1(F) = -r_1(F) = -\frac12 B_2(3/2) +\frac12 \Big(H_F(0) + 2H_F(1) +2H_F(2)\Big) = \frac{13}{24} +\xi_F + H_F(2),
\]
and the result follows in view of \eqref{1-1}. \qed

\medskip
\section{Virtual $\gamma$-factors} 

\smallskip
With the cases of holomorphic and non-holomorphic modular forms of level 1 in mind, we define the {\it virtual $\gamma$-factors} as
\begin{equation}
\label{5-1}
\gamma(s) = 
\begin{cases}
(2\pi)^{-s} \Gamma(s+\mu) &\text{with $\mu>0$} \\
\pi^{-s} \Gamma\big(\frac{s+\epsilon+i\kappa}{2}\big) \Gamma\big(\frac{s+\epsilon-i\kappa}{2}\big) &\text{with $\epsilon\in\{0,1\}$ and $\kappa\geq 0$.}
\end{cases}
\end{equation}
Obviously, we say that the first $\gamma$-factor in \eqref{5-1} is of {\it Hecke type} and the second is of {\it Maass type}.
In view of Remark 2.1, the virtual $\gamma$-factors have degree 2 and conductor 1. Moreover, by \eqref{2-5} we also have
\begin{equation}
\label{5-2}
\xi_\gamma= 
\begin{cases}
2\mu-1 \\
2(\ep-1)
\end{cases}
\text{and} \quad\chi_\gamma:=  \xi_\gamma + H_\gamma(2) +\frac23= 
\begin{cases}
2\mu^2 \\
- 2\kappa^2,
\end{cases}
\end{equation}
and clearly, with obvious notation,
\begin{equation}
\label{5-3}
\prod_{j=1}^r \lambda_j^{2i\Im{\mu_j}} =1
\end{equation}
for all virtual $\gamma$-factors.

\smallskip
As suggested by the name, not every virtual $\gamma$-factor corresponds to an existing $L$-function. The list of all normalized $F\in\S^\sharp$ whose $\ga$-factor is a virtual $\gamma$-factor is provided by the classical converse theorems of Hecke \cite{Hec/1983} and Maass \cite{Maa/1949}. In the next lemma, ``eigenvalue of the Laplacian'' means eigenvalue of the hyperbolic Laplacian on the real analytic, $L^2$ and $\Gamma_0(1)$-invariant functions on the upper half-plane $\HH$.

\medskip
{\bf Lemma 5.1.} {\sl If the $\gamma$-factor of $F\in\S^\sharp$ is a virtual $\gamma$-factor and $F$ is normalized, then one of the following cases holds true:

(i) $\mu= \frac{k-1}{2}$ with an even integer $k\geq 12$, in which case $\omega_F=(-1)^{k/2}$ and $F(s) = L(s+\mu,f)$ with some holomorphic cusp form $f$ of level $1$ and weight $k$;

(ii) $\rho=1/4 + \kappa^2$ is an eigenvalue of the Laplacian, in which case $\omega_F = (-1)^\ep$ and $F(s) = L(s,u)$ with some Maass form $u$ of level $1$, weight $0$, parity $\ep$ and with eigenvalue $\rho$.} 

\medskip
{\it Proof.} We deal first with virtual $\gamma$-factors of Hecke type. Suppose that a normalized $F\in\S^\sharp$ satisfies
\begin{equation}
\label{5-4}
\gamma(s)F(s)=\omega \gamma(1-s)F(1-s).
\end{equation}
Then $a(n) \ll n^c$ for some $c>0$ and, by the change of variable $s\mapsto s-\mu$, \eqref{5-4} becomes 
\[
(2\pi)^{-s} \Gamma(s)G(s) = \omega (2\pi)^{-(k-s)} \Gamma(k-s)G(k-s),
\]
where
\[
k=2\mu+1>1 \quad \text{and} \quad G(s) = F(s-\mu).
\]
This is Hecke's functional equation with signature $(1,k,\omega)$, see Definition 2.1 of Berndt-Knopp \cite{Be-Kn/2008}. Moreover, by \eqref{5-3} we have
\begin{equation}
\label{5-5}
\omega_F=\omega.
\end{equation}
But, by Lemma 4.1 and \eqref{5-2}, $\xi_F=2\mu-1$ is an even integer, thus $k$ is also an even integer and $\omega=(-1)^{k/2}$ by \eqref{5-5}. Hence $G(s)$ has signature
\[
(1,k,(-1)^{k/2}) \quad \text{with an even integer $k\geq2$},
\]
therefore by Hecke's converse theorem, see Theorem 2.1 of \cite{Be-Kn/2008}, $G(s) = L(s,f)$, the $L$-function associated with a modular form $f$ of level 1 and weight $k$. According to Theorem 4 of Chapter VII of Serre \cite{Ser/1973} (where the weight of $f$ is denoted by $2k$) the space of modular forms of weight $k\in2\ZZ$ is trivial for $k<0$ and $k=2$, while it has dimension 1 for $k=0$, 4, 6, 8, 10 and is generated, respectively, by 1 and the Eisenstein series $G_4$, $G_6$, $G_8$, $G_{10}$. Moreover, for $k\geq12$ the subspace of cusp forms has always dimension $\geq 1$. But $L(s,G_k)$ has a pole at $s=k$, thus $F$ should have a pole at $s=(k+1)/2\neq1$, impossible for $k\geq4$. On the other hand, if $f$ is a Hecke eigenform of weight $k\geq12$ then $F(s)=L(s+\frac{k-1}{2},f)$ is normalized, and our first assertion follows.

\smallskip
Next we deal with virtual $\gamma$-factors of Maass type. The argument is similar, this time based on the version allowing poles of Maass' converse theorem given by Theorem 2.1 of Raghunathan \cite {Rag/2010}, with $\nu=1/2+i\kappa$. Suppose that a normalized $F$ satisfies \eqref{5-4} with such a $\gamma$-factor. Then by \eqref{5-5}, Lemma 4.1 and \eqref{5-2} we have $\omega = (-1)^\ep$. Hence our assertion follows from the above quoted version of Maass' converse theorem, since $\nu(1-\nu)=1/4+\kappa^2=\rho$ must be an eigenvalue of the Laplacian, thus its eigenspace is nontrivial and therefore we may choose a Maass form $u$ associated with $\rho$ and with parity $\ep$ such that $F(s)=L(s,u)$ is normalized. \qed

\medskip
{\bf Remark 5.1.} The limitations $\mu>0$ and $\kappa\geq0$ imposed on the virtual $\gamma$-factors in \eqref{5-1} are actually natural. Indeed, we may restrict to $\kappa\geq0$ since there are no exceptional eigenvalues of the Laplacian in the case of $\Gamma_0(1)$, see Selberg \cite{Sel/1965}. Moreover, if $\mu\leq0$ then the weight $k$ in the above proof is an even integer $\leq 0$, thus by the above quoted Theorem 4 of \cite{Ser/1973} the space of modular forms is either trivial (if $k<0$) or generated by $1$ (if $k=0$), in which case $L(s,1)$ is identically vanishing and hence $F(s)\equiv0$, impossible. \qed

\medskip
By \eqref{2-8}, given any virtual $\gamma$-factor in \eqref{5-1} we define the {\it virtual $h$-function} as
\begin{equation}
\label{5-6}
h_\gamma(s) =
\begin{cases}
 (2\pi)^{2(s-1)} \Gamma(1-s+\mu) \Gamma(1-s-\mu) \\
 \frac14 \pi^{2s-3} \Gamma\Big(\frac{1-s+\epsilon-i\kappa}{2}\Big)  \Gamma\Big(\frac{2-s-\epsilon+i\kappa}{2}\Big) \Gamma\Big(\frac{1-s+\epsilon+i\kappa}{2}\Big) \Gamma\Big(\frac{2-s-\epsilon-i\kappa}{2}\Big).
\end{cases}
\end{equation}
Hence by the reflection formula for the $\Gamma$ function we have that
\begin{equation}
\label{5-7}
h_\gamma(s) S_\gamma(s) = \frac{\gamma(1-s)}{\gamma(s)},
\end{equation}
and by \eqref{2-12} and \eqref{2-9}
\begin{equation}
\label{5-8}
h_\gamma(s) \approx \frac{(4\pi)^{2s-1}}{\sqrt{2\pi}} \sum_{\ell=0}^\infty d_\ell(\gamma) \Gamma(3/2 -2s-\ell), \quad d_0(\gamma)=1.
\end{equation}

\medskip
{\bf Remark 5.2.} Since $\ep\in\{0,1\}$, the virtual $h$-functions do not depend on $\ep$ thanks to the symmetry of the $\Gamma$-factors in \eqref{5-6}. In particular, this shows that the $h$-functions do not determine uniquely the $\gamma$-factors. \qed

\medskip
Although virtual $\gamma$-factors do not correspond, in general, to some $F\in\S^\sharp$, their structural invariants $d_\ell(\gamma)$ satisfy the same properties of the invariants $d_\ell(F)$ of the normalized functions $F\in\S^\sharp$ of degree 2 and conductor 1.

\medskip
{\bf Lemma 5.2.} {\sl The structural invariants $d_\ell(\gamma)$ of virtual $\gamma$-factors satisfy the same properties of the invariants $d_\ell(F)$ in Proposition $4.1$, Corollary $4.1$ and Lemma $4.3$.}

\medskip
{\it Proof.} We first show that there exist polynomials $P_\ell,Q_\ell\in\RR[x]$ such that
\begin{equation}
\label{5-9}
d_\ell(\gamma) = 
\begin{cases}
P_\ell(\mu) \\
Q_\ell(\kappa)
\end{cases}
\end{equation}
for every virtual $\gamma$-factor in \eqref{5-1}. Note, in view of Remark 5.2, that the polynomials $Q_\ell$, if they exist, do not depend on $\ep$. Now we observe that Lemma 4.3 is proved by Stirling's formula, without using the properties of $F\in\S^\sharp$. Hence in the present case \eqref{4-39} becomes
\begin{equation}
\label{5-10}
d_\ell(\gamma) =  \sum_{h=1}^\ell A_{h,\ell} V_h(\gamma),
\end{equation}
where $A_{h,\ell}$ is defined by \eqref{4-38} and $V_h(\gamma)$ by \eqref{4-37} and \eqref{4-35}, clearly with $\gamma$ in place of $F$. But the $H$-invariants in \eqref{4-35} are defined in terms of Bernoulli polynomials involving $\mu$ and $\kappa$, hence the dependence of $d_\ell(\gamma)$ on $\mu$ and $\kappa$ is polynomial as well. Moreover, $A_{h,\ell}\in\RR$ so $P_\ell\in\RR[x]$, and $Q_\ell(\kappa)$ involves Bernoulli polynomials at $i\kappa$ and $-i\kappa$, thus its coefficients are also real. Hence \eqref{5-9} is proved.

\smallskip
Next we observe that if a virtual $\gamma$-factor is actually a $\gamma$-factor of a function $F\in\S^\sharp$ as in Theorem 1.1, then by Proposition 4.1 we have
\begin{equation}
\label{5-11}
Q_N(d_0(\gamma),\dots,d_N(\gamma))=0.
\end{equation}
But $Q_N(d_0(\gamma),\dots,d_N(\gamma))$ is a polynomial in $\mu$ or $\kappa$ by \eqref{5-9}, and by Lemma 5.1 there exist infinitely many $\mu$ and $\kappa$ for which \eqref{5-11} holds (remember that the spectrum of the Laplacian is infinite). Thus \eqref{5-11} holds identically in $\mu$ and $\kappa$, i.e. Proposition 4.1 holds for all virtual $\gamma$-factors. Hence the analog of Corollary 4.1 holds as well, and the analog of Lemma 4.3 follows from \eqref{5-10}. \qed

\medskip
Now we are ready to prove the main result of this section.

\medskip
{\bf Proposition 5.1.} {\sl Let $F\in\S^\sharp$ be as in Theorem $1.1$. Then there exists a virtual $\gamma$-factor such that
\[
\chi_F= \chi_\gamma \quad \text{and} \quad h_F(s) = \omega_F h_\gamma(s).
\]
Moreover, such a $\gamma$-factor is uniquely determined if we choose $\ep$ in \eqref{5-1} satisfying $\omega_F=(-1)^\ep$, and it is of Hecke or Maass type depending on $\chi_F>0$ or $\chi_F\leq 0$, respectively.}

\medskip
{\it Proof.} Thanks to Lemmas 4.3 and 5.2 we have
\begin{equation}
\label{5-12}
d_1(F) = d_1(\gamma) \ \Longleftrightarrow \ \chi_F= \chi_\gamma.
\end{equation}
Moreover, by Lemmas 4.1 and 4.3 we have $d_1(F)\in\RR$ (recall that $\xi_F=H_F(1)$), and by \eqref{5-2} the set of values of $\chi_\gamma$ coincides with $\RR$. Thus there exists a virtual $\gamma$-factor such that $d_1(F) = d_1(\gamma)$. But then, thanks to Corollary 4.1 and Lemma 5.2, we have that
\[
d_\ell(F) = d_\ell(\gamma) \quad \text{for every $\ell\geq0$,}
\]
and the first assertion of the proposition follows in view of \eqref{2-12} and \eqref{5-8}. Indeed, if $d_\ell(F) = d_\ell(\gamma)$ for every $\ell\geq0$ then the argument after Lemma 2.2. in \cite{Ka-Pe/2002}, applied to $h_F(s)$ and $h_\ga(s)$ in place of $\ga_1(s)$ and $\ga_2(s)$, shows that $h_F(s) = \omega_F h_\gamma(s)$. For a proof of Lemma 2.2 in \cite{Ka-Pe/2002} we refer to Lemma D in \cite{Ka-Pe/revisited}.

\smallskip
In view of \eqref{5-2} and the equality in the right hand side of \eqref{5-12}, if $\chi_F>0$ there is only one choice for the virtual $\gamma$-factor of Hecke type. Instead, if $\chi_F\leq 0$ there are two possible choices of virtual $\gamma$-factors of Maass type, depending on the value of $\ep$. Thus, since $\omega_F=\pm1$ by Lemma 4.1, we choose $\ep\in\{0,1\}$ in \eqref{5-1} in such a way that $\omega_F=(-1)^\ep$, and the second assertion follows as well. \qed

\medskip
Recalling the definition of the functions $S_F(s)$ and $S_\gamma(s)$ in Section 2.1 and writing
\begin{equation}
\label{5-13}
R(s) = R_{F,\ga}(s) := \frac{S_F(s)}{S_\gamma(s)}
\end{equation}
we have the following important consequence of Proposition 5.1.

\medskip
{\bf Corollary 5.1.}  {\sl Let $F\in\S^\sharp$ be as in Theorem $1.1$ and $R(s)$ be as in \eqref{5-13}. Then there exists a unique virtual $\gamma$-factor, in the sense of Proposition $5.1$, such that $F$ satisfies the functional equation}
\[
\gamma(s) F(s) = \omega_F R(s) \gamma(1-s) F(1-s).
\]

\medskip
{\it Proof.} We start with the functional equation of $F$ written as in \eqref{2-10} and then use Lemma 4.1, Proposition 5.1 and \eqref{5-7} to get
\[
F(s) = \omega_F h_\gamma(s) S_F(s) F(1-s) = \omega_F \frac{\gamma(1-s)}{\gamma(s)} \frac{S_F(s)}{S_\gamma(s)} F(1-s);
\]
thus the result follows. \qed

\medskip
We conclude the section with a study of the function $R(s)$ in \eqref{5-13}. We shall always assume that $\gamma$ in \eqref{5-13} is the virtual $\gamma$-factor associated with $F$ in the sense of Proposition 5.1. Moreover, with a convenient abuse of notation, we rewrite \eqref{2-7} as
\begin{equation}
\label{5-14}
S_F(s) = \sum_{j=0}^N a_j e^{i\pi\omega_js}
\end{equation}
with
\begin{equation}
\label{5-15}
N\geq 1, \quad -1=\omega_0<\omega_1<\cdots<\omega_N=1, \quad \omega_j=-\omega_{N-j}, \quad a_j\neq0.
\end{equation}
Further, a computation based on \eqref{2-7} shows that
\begin{equation}
\label{5-16}
a_0=a_N=-\omega_F,
\end{equation}
since $\omega_F=-e^{i\pi\xi_F/2}$ and $\xi_F\in2\ZZ$ by Lemma 4.1.

\medskip
{\bf Lemma 5.3.} {\sl Let $F(s)$, $R(s)$ and $S_F(s)$ be as in Theorem $1.1$, \eqref{5-13} and \eqref{5-14}, respectively.

i) If $R(s)$ is not constant, then $N\geq3$ and $\omega_{N-1}>0$.

ii) If $\ga$ is of Maass type we have
\[
R(s) = 1 + O\big(e^{-\delta|t|}\big)
\]
as $|t|\to\infty$, with some $\delta>0$.}

\medskip
{\it Proof.} i) We first note that applying twice the functional equation in Corollary 5.1 we get
\begin{equation}
\label{5-17}
R(s)R(1-s)=1.
\end{equation}
Moreover, for $N=1$ and $N=2$ we have
\begin{equation}
\label{5-18}
S_F(s) = -\omega_F \big(e^{\pi is} + e^{-\pi is}\big) +a = -2\omega_F\cos(\pi s) + a
\end{equation}
with $a=0$ if $N=1$ and $a=a_1\neq0$ if $N=2$. This is immediate from \eqref{5-14},\eqref{5-15} and \eqref{5-16} observing, if $N=2$, that $\omega_1=0$ since $\omega_1=-\omega_1$.

\smallskip
We start with the case of Hecke type virtual $\gamma$-factors, where
\begin{equation}
\label{5-19}
S_\ga(s) = e^{\pi i(\mu-1/2)} e^{\pi is} + e^{-\pi i(\mu-1/2)} e^{-\pi is}.
\end{equation}
From \eqref{5-17},\eqref{5-18} and \eqref{5-19} we obtain
\begin{equation}
\label{5-20}
\frac{S_F(s)S_F(1-s)}{S_\ga(s)S_\ga(1-s)} = \frac{-\omega_F \big(e^{\pi is} + e^{-\pi is}\big) +a}{e^{\pi i(\mu-1/2)} e^{\pi is} + e^{-\pi i(\mu-1/2)} e^{-\pi is}} 
\frac{-\omega_F \big(e^{\pi is} + e^{-\pi is}\big) -a}{e^{\pi i(\mu-1/2)} e^{-\pi is} + e^{-\pi i(\mu-1/2)} e^{\pi is}} =1
\end{equation}
which implies, since $\omega_F=\pm1$ thanks to Lemma 4.1, that
\begin{equation}
\label{5-21}
2-a^2 = 2\cos\big(2\pi(\mu-1/2)\big).
\end{equation}
If $N=1$ then \eqref{5-21} implies that $\mu-1/2\in\ZZ$ and hence \eqref{5-19} becomes
\[
S_\ga(s) = \pm1 \big(e^{\pi is} + e^{-\pi is}\big),
\]
thus $R(s)$ is constant in this case. Let now $N=2$ and suppose that $S_\ga(s)=0$. Then by \eqref{5-19} we have that $t=0$ and $\si= -\mu+k$ for some $k\in\ZZ$, hence from \eqref{5-20} we deduce that
\begin{equation}
\label{5-22}
\cos^2(\pi\mu) = a^2.
\end{equation}
Since $\cos(x-\pi) = -\cos x$ and $\cos(2x) = 2\cos^2x-1$, equation \eqref{5-21} can be written as $a^2=4\cos^2(\pi\mu)$, which contradicts \eqref{5-22} since $a\neq0$. Thus we cannot have $N=2$, and the result follows in the Hecke case.

\smallskip
For virtual $\gamma$-factors of Maass type we have
\begin{equation}
\label{5-23}
S_\gamma(s) = 4 \sin\big(\pi(\frac{s+\ep+i\kappa}{2})\big) \sin\big(\pi(\frac{s+\ep-i\kappa}{2})\big) = -2(-1)^\ep \cos(\pi s) + a(\gamma)
\end{equation}
with a certain $a(\gamma)$. If $N=1$, from \eqref{5-17},\eqref{5-18} and the first equation in  \eqref{5-23} we get
\[
\begin{split}
16 \sin\big(\pi(\frac{s+\ep+i\kappa}{2})\big) \sin\big(\pi(\frac{s+\ep-i\kappa}{2})\big) &\sin\big(\pi(\frac{1-s+\ep+i\kappa}{2})\big) \sin\big(\pi(\frac{1-s+\ep-i\kappa}{2})\big) \\
&= 4\cos^2(\pi s),
\end{split}
\]
impossible since zeros on the two sides do not match. Thus we cannot have $N=1$. If $N=2$, then by \eqref{5-17},\eqref{5-18} and the second expression for $S_\gamma(s)$ in \eqref{5-23} we have
\[
\begin{split}
\big(-2\omega_F\cos(\pi s) + a\big) &\big(-2\omega_F\cos(\pi (1-s)) + a\big) \\
&= \big(-2(-1)^\ep \cos(\pi s) + a(\gamma)\big)(-2(-1)^\ep \cos(\pi (1-s)) + a(\gamma)\big),
\end{split}
\]
which gives $a^2 = a(\gamma)^2$ and hence $a=\pm a(\gamma)$. Recalling that $\omega_F=(-1)^\ep$ by Proposition 5.1, if $a=a(\ga)$ we have $S_F(s)=S_\gamma(s)$ by \eqref{5-18} and \eqref{5-23}, hence $R(s)$ is constant in this case. If $a=-a(\ga)$, again from \eqref{5-18} and \eqref{5-23} we have $S_\gamma(s) = -S_F(1-s)$ and hence Corollary 5.1 yelds
\[
S_F(1-s) \gamma(s) F(s) =- \omega_F S_F(s) \gamma(1-s) F(1-s).
\]
Dividing this by the functional equation \eqref{2-1} of $F$, denoting by $\gamma_F(s)$ the $\gamma$-factor of $F$ and recalling Lemma 4.1 we obtain
\[
\frac{S_F(1-s)\gamma(s)}{\gamma_F(s)} = - \frac{S_F(s)\gamma(1-s)}{\gamma_F(1-s)}.
\]
Thus the function $f(s) = S_F(1-s) \gamma(s) \gamma_F(1-s)$ satisfies $f(s)=-f(1-s)$, and in particular $f(1/2)=0$. Therefore 
\[
S_F(1/2) \gamma(1/2) \gamma_F(1/2)=0,
\]
and since $\gamma$-factors do not vanish this implies $0=S_F(1/2)=a$, a contradiction. The proof of i) is now complete. 

\smallskip
ii) A computation based on \eqref{2-7} shows that as $|t|\to\infty$
\[
S_F(s) = e^{-\text{sgn}(t)\pi i(s+\xi_F/2)} \big(1+ O(e^{-c_1|t|})\big)
\]
with some $c_1>0$. Hence by Lemma 4.1 this becomes
\[
S_F(s) =-\omega_F  e^{-\text{sgn}(t)\pi is} \big(1+ O(e^{-c_1|t|})\big)
\]
since $\omega_F=\pm1$. In view of \eqref{5-23}, for virtual $\gamma$-factors of Maass type we have
\[
S_\gamma(s) = -(-1)^\ep  e^{-\text{sgn}(t)\pi is} \big(1+ O(e^{-c_2|t|})\big)
\]
with some $c_2>0$. Hence ii) follows, since in this case $\omega_F=(-1)^\ep$ by Proposition 5.1. \qed

\medskip
\section{Period functions}

\smallskip
For a function $F$ as in Theorem $1.1$ and $z$ in the upper half-plane $\HH$ we define 
\begin{equation}
\label{6-1}
f(z) = \sum_{n=1}^\infty a(n) n^\lambda e(nz), \quad \text{where} \quad \lambda=
\begin{cases}
\mu & \text{if $\chi_F >0$} \\
i\kappa & \text{if $\chi_F \leq 0$,}
\end{cases}
\end{equation}
$\mu,\kappa$ are as in \eqref{5-1} and $\chi_F$ is defined by \eqref{1-1}. Clearly, $f(z)$ is holomorphic on $\HH$. We shall repeatedly use the following simple, and essentially well known, lemma.

\medskip
{\bf Lemma 6.1.} {\sl For every $F\in\S^\sharp$ with positive degree and $\la\in\CC$ the function $f(z)$ cannot be continued to an entire function.}

\medskip
{\it Proof.} From the convergence properties of the Dirichlet series of $F$ we have that $a(n)\ll n^{1+\ep}$, hence the power series $g(z) = \sum_{n=1}^\infty a(n)n^\la z^n$ is convergent for $|z|<1$. If $f(z)$ is entire, then $g(z)$ is also entire since clearly it is bounded around $z=0$. Therefore $g(z)$ is convergent for some $|z_0|>1$, thus $a(n)n^\la \ll |z_0|^{-n}$ and hence $a(n)\ll n^{-A}$ for every $A>0$. As a consequence the Dirichlet series of $F$ is everywhere convergent, a contradiction since its Lindel\"of function $\mu_F(\si)$ is positive for $\si<0$. \qed

\medskip
In the next two propositions we study the function $f(z)$, and other related functions, distinguishing the cases $\chi_F>0$ and $\chi_F\leq0$, corresponding to associated virtual $\gamma$-factors of Hecke and Maass type, since the arguments are somewhat different. Actually, the Hecke case is definitely more delicate. In what follows the suffixes $H$ and $M$ are related to the Hecke or Maass cases, respectively.

\smallskip
Let $R(s)$ be as in \eqref{5-13}, $\gamma$ being the virtual $\gamma$-factor associated with $F$, and
\begin{equation}
\label{6-2}
\begin{split}
0<c_1<1 \ \text{if} \ \la=i\kappa \ & \text{and} \ 0<c_1<\mu \ \text{if} \ \la=\mu,  \\
 c_2=1+2\Re(\la)-c_1 \ & \text{and} \ \delta = \mu-c_1;
 \end{split}
\end{equation}
clearly, $\delta>0$. We define the following functions.

(a) If $\la=i\kappa$, i.e. $\gamma$ is of Maass type, we write
\begin{equation}
\label{6-3}
P_M(z) = \frac{i^{\ep+1}\cos(\pi\lambda)}{2\pi^{\lambda+3/2}} \frac{1}{2\pi i} \int_{(c_2)} \Gamma\big(\frac{1-\ep+s}{2}\big) \Gamma\big(\frac{2-\ep-s+2\lambda}{2}\big) \ga(s-\lambda) F(s-\lambda) z^{s-1-2\lambda} \d s
\end{equation}
and
\begin{equation}
\label{6-4}
Q_M(z) = \frac{\omega_F}{2\pi i} \int_{(c_1)} \big(R(s-\lambda)-1\big) \Gamma(s) \frac{\ga(1-s+\lambda)}{\ga(s-\lambda)} F(1-s+\lambda) (-2\pi iz)^{-s} \d s.
\end{equation}
Clearly, the integral in \eqref{6-3} is absolutely convergent for $|\arg(z)|<\pi$ and $P(z)$ is holomorphic there, while thanks to Lemma 5.3 the integral in \eqref{6-4} is certainly absolutely convergent for $z\in\HH$ and $Q_M(z)$ is holomorphic there.

(b) If $\la=\mu$, i.e. $\ga$ is of Hecke type, we write
\begin{equation}
\label{6-5}
\begin{split}
Q_H(z) = \omega_F(2\pi)^{1-\mu} (-iz)^{-\mu-1} \sum_{n=1}^\infty &a(n) \frac{1}{2\pi i} \int_{1+\delta} \Big(\frac{1}{S_F(s)}+\frac{1}{2\omega_F\cos(\pi s)}\Big) \\
&\times \frac{1}{\Ga(1-s-\mu)} \Big(\frac{2\pi in}{z}\Big)^{-s} \d s.
\end{split}
\end{equation}
The analytic properties of $Q_H(z)$ are established in Proposition 6.1, but it is not difficult to see that $Q_H(z)$ is at least holomorphic for $z\in\mathbb H$ using \eqref{6-12} below.

(c) For every $\la$ as above
\begin{equation}
\label{6-6}
L(z) =  \res_{s=1+\lambda} \, \Gamma(s) F(s-\lambda)(-2\pi iz)^{-s};
\end{equation}
clearly, $L(z)\equiv0$ if $F$ is entire and otherwise is holomorphic for $|\arg(z)|<\pi$, since it is a linear combinations of terms of type $z^a$ and $z^a\log^kz$ with $a\in\CC$ and $k\in\NN$. Finally, for $z\in\HH$ we write
\begin{equation}
\label{6-7}
\psi(z) = f(z) -z^{-2\la-1}f(-1/z).
\end{equation}

\medskip
{\bf Proposition 6.1.} {\sl Let $F$ be as in Theorem $1.1$ with $\chi_F>0$. Then with the above notation
\[
\psi(z) = Q_H(z) + P_H(z) + L(z)
\]
and $Q_H(z)$ is holomorphic for $|\arg(z)|<\pi$, where $P_H(z)$ is a certain holomorphic function for $|\arg(z)|<\pi$.}

\medskip
{\it Proof.} For future reference we take the first steps of the proof keeping $\chi_F$ and $\la$ general. Let $z=x+iy\in\HH$. Since $F(s-\lambda)$ is absolutely convergent for $\si>1+\Re(\lambda)$, by Mellin's transform we have
\begin{equation}
\label{6-8}
f(iy) =  \sum_{n=1}^\infty a(n) n^\lambda e^{-2\pi ny} = \frac{1}{2\pi i} \int_{(c)} (2\pi)^{-s} \Gamma(s) F(s-\lambda) y^{-s} \d s,
\end{equation}
where $c> \Re(\lambda)+1$.  Recalling \eqref{6-6} and that $F$ has at most a pole at $s=1$, we shift the integration line to $\si= c_1$ as in \eqref{6-2}, thus getting
\[
f(iy) = \frac{1}{2\pi i} \int_{(c_1)} (2\pi)^{-s} \Gamma(s) F(s-\lambda) y^{-s} \d s + L(iy).
\]
Hence by Corollary 5.1 we obtain
\begin{equation}
\label{6-9}
f(iy) =  \frac{\omega_F}{2\pi i} \int_{(c_1)} (2\pi)^{-s} \Gamma(s) R(s-\lambda) \frac{\ga(1-s+\lambda)}{\ga(s-\lambda)} F(1-s+\lambda) y^{-s} \d s + L(iy).
\end{equation}

\smallskip
Now we recall that if $\chi_F>0$ then $\la=\mu>0$ and the associated virtual $\gamma$-factor is $\gamma(s) = (2\pi)^{-s}\Gamma(s+\mu)$; in what follows we restrict ourselves to this case. Thus, for $\chi_F>0$, by the change of variable $s\mapsto 1-s+\mu$ we see that \eqref{6-9} becomes
\begin{equation}
\label{6-10}
f(iy) = \omega_F  (2\pi)^{-\mu} y^{-\mu-1} \frac{1}{2\pi i} \int_{(1-c_1+\mu)} R(1-s) \gamma(s) F(s) y^{s} \d s + L(iy).
\end{equation}
But by \eqref{6-2} we have that $1-c_1+\mu=1+\delta$; for later use, we may assume that $\delta>0$ is sufficiently small. Moreover, thanks to \eqref{5-13} and \eqref{5-17} we rewrite \eqref{6-10} as
\[
f(iy) =  \omega_F  (2\pi)^{-\mu} y^{-\mu-1} \frac{1}{2\pi i} \int_{(1+\delta)} \frac{S_\ga(s)}{S_F(s)} \gamma(s) F(s) y^{s} \d s + L(iy).
\]
Hence, expanding $F(s)$, recalling that $S_\ga(s) = 2\sin(\pi(s+\mu))$ by definition and applying the reflection formula to $\ga(s)$, we finally obtain that
\begin{equation}
\label{6-11}
f(iy) = \omega_F (2\pi)^{1-\mu} y^{-\mu-1} \sum_{n=1}^\infty a(n) \frac{1}{2\pi i} \int_{(1+\delta)} \frac{1}{S_F(s)\Ga(1-s-\mu)}  \Big(\frac{2\pi n}{y}\Big)^{-s} \d s + L(iy).
\end{equation}

\smallskip
From \eqref{5-14}-\eqref{5-16} we have that
\[
S_F(s) = -2\omega_F \cos(\pi s) + \sum_{j=1}^{N-1} a_j e^{i\pi\omega_js}
\]
with $|\omega_j|<1$ for $1\leq j\leq N-1$. Thus, as $|t|\to\infty$
\begin{equation}
\label{6-12}
\frac{1}{S_F(s)} +\frac{1}{2\omega_F\cos(\pi s)} \ll e^{-\pi(1+\rho)|t|}
\end{equation}
for some $0<\rho<1/2$. Hence we rewrite \eqref{6-11} as
\begin{equation}
\label{6-13}
f(iy) = -\frac12 (2\pi)^{1-\mu} y^{-\mu-1} \sum_{n=1}^\infty a(n) J\Big(\frac{2\pi n}{y}\Big) + H(y) + L(iy)
\end{equation}
where, introducing a new variable $w$ for later use,
\[
J(w) = \frac{1}{2\pi i} \int_{(1+\delta)} \frac{1}{\cos(\pi s)\Gamma(1-s-\mu)} w^{-s}\d s
\]
and
\[
H(y) =  \omega_F (2\pi)^{1-\mu} y^{-\mu-1} \sum_{n=1}^\infty \frac{a(n)}{2\pi i} \int_{(1+\delta)} \Big\{\frac{1}{S_F(s)} + \frac{1}{2\omega_F\cos(\pi s)}\Big\} \frac{1}{\Ga(1-s-\mu)} \Big(\frac{2\pi n}{y}\Big)^{-s} \d s.
\]

\smallskip
Now we consider $y$ as a complex variable. Since $z=x+iy$ and our aim is to obtain analytic continuation to $|\arg(z)|<\pi$, we are mainly interested in the range $y\in \CC\setminus [0,i\infty)$ with $-3\pi/2<\arg(y)<\pi/2$. From \eqref{6-12} and Stirling's formula we see that the sum of integrals defining $H(y)$ is holomorphic for $|\arg(y)| < \pi/2+\pi \rho$. Moreover, by Stirling's formula the integral defining $J(2\pi n/y)$ is holomorphic for $|\arg(y)|<\pi/2$, and clearly
\[
J(2\pi n/y) \ll_y n^{-1-\delta}.
\]
Thus the series on the right hand side of \eqref{6-13} is also holomorphic for $|\arg(y)|<\pi/2$.

\smallskip
The next step is a closer study of $J(w)$. Since later on we choose $w = 2\pi n/y$, we consider the domain $\D = \{w\in \CC\setminus [0,-i\infty)$ with $-\pi/2< \arg(w)<3\pi/2\}$. We shift the line of integration in $J(w)$ to $-\infty$ and, recalling that $\delta>0$ is sufficiently small, we collect the residues of the simple poles at the points $s=-\ell+1/2$ with $\ell\geq0$. Thus we obtain that
\begin{equation}
\label{6-14}
J(w) = -\frac{w^{-1/2}}{\pi} \sum_{\ell=0}^\infty \frac{(-w)^\ell}{\Gamma(\ell + 1/2 -\mu)} = -\frac{w^{-1/2}}{\pi} E_{1/2-\mu}(w),
\end{equation}
say, where the branch of $\log w$ in $w^{-1/2}$ is chosen as above. Formula \eqref{6-14} holds for $|\arg(w)|<\pi/2$, since by Stirling's formula the integral tends to 0 as the integration line is shifted to $-\infty$. Moreover, the series in \eqref{6-14} is everywhere convergent, thus $E_{1/2-\mu}(w)$ is an entire function and hence \eqref{6-14} provides the holomorphic continuation of $J(w)$ to the domain $\D$. Actually, $E_{1/2-\mu}(w)$ is a special case of a {\it two-parametric Mittag-Leffler function}, see Chapter 4 of Gorenflo-Kilbas-Mainardi-Rogosin \cite{GKMR/2014}. More precisely we have that
\[
E_{\beta}(w) = E_{1,\beta}(-w),
\]
where $\beta=1/2-\mu$ and $E_{1,\beta}(w)$ is the function defined by (4.1.1) of \cite{GKMR/2014}. One easily checks (see also (4.3.2) of \cite{GKMR/2014} with $m=n=1$) that $E_\beta(w)$ satisfies the first order linear differential equation
\[
E_\beta'(w) = -\Big(1+\frac{\beta-1}{w}\Big) E_\beta(w) + \frac{1}{w\Gamma(\beta-1)}.
\]
Writing $w=u+iv$, we solve such a differential equation by the standard method, choosing an initial point $w_0\neq 0$ of the form $u_0+iv$ and then letting $u_0\to-\infty$. Thus we obtain for $w\in\D$ that
\begin{equation}
\label{6-15}
E_\beta(w) = \kappa_0 e^{-w} w^{1-\beta} + \frac{ e^{-w} w^{1-\beta} }{\Ga(\beta-1)} \int_{\LL(w)} e^z z^{\beta-2} \d z
\end{equation}
with some constant $\kappa_0$, where $\LL(w)$ is a path from $-\infty$ to $w$ not crossing the half-line $[0,-i\infty)$. Precisely, if $0<\arg(w)<3\pi/2$ then we choose $\LL(w)$ as $(-\infty+i\Im(w),w]$, while if $-\pi/2< \arg(w)\leq 0$ then $\LL(w)$ is chosen as $(-\infty+i, \Re(w)+i] \cup [\Re(w)+i,w]$; note that the integral in \eqref{6-15} is absolutely convergent. 

\smallskip
We remark that formula \eqref{6-15} follows also from (4.3.6) of \cite{GKMR/2014} with $z=-w$, $z_0=-w_0$ and $n=1$, after the change of variable $\tau=-\xi$ inside the integral. Indeed, the constant term in (4.3.6) is convergent as $u_0\to-\infty$ thanks to (4.4.16) in Theorem 4.3 of \cite{GKMR/2014}.

\smallskip
In view of \eqref{6-1}, from \eqref{6-13}-\eqref{6-15} with $w=2\pi n/y$ we get
\begin{equation}
\label{6-16}
f(iy) = \kappa_1 (iy)^{-2\mu-1} f(i/y) + \frac{y^{-2\mu-1}}{\Ga(-1/2-\mu)} \widetilde{f}(y) + H(y) + L(iy),
\end{equation}
where $\kappa_1= \kappa_0i^{2\mu+1}$ and
\begin{equation}
\label{6-17}
 \widetilde{f}(y) = \sum_{n=1}^\infty a(n) n^\mu e^{-2\pi n/y} \int_{\LL(2\pi n/y)} e^z z^{-\mu-3/2} \d z.
\end{equation}
Note that the integral in \eqref{6-17} is
\[
\ll e^{2\pi n\Re(y)/|y|^2} n^{-\mu-3/2}
\]
uniformly for $y$ in any compact subset of $\CC\setminus [0,i\infty)$, therefore the series in \eqref{6-17} is absolutely convergent for the same values of $y$ and hence $\widetilde{f}(y)$ is holomorphic for $y\in\CC\setminus [0,i\infty)$. 

\smallskip
Coming back to the variable $z=x+iy\in\HH$, in view of \eqref{6-5} we rewrite \eqref{6-16} as
\begin{equation}
\label{6-18}
\begin{split}
f(z) &= \kappa_1 z^{-2\mu-1} f(-1/z) + \frac{(-iz)^{-2\mu-1}}{\Ga(-1/2-\mu)} \widetilde{f}(-iz) + Q_H(z)   +L(z)  \\
&=\kappa_1 z^{-2\mu-1} f(-1/z) + \widetilde{\psi}(z),
\end{split}
\end{equation}
say, where $\widetilde{\psi}(z)$ is holomorphic for $-\pi\rho < \arg(z)<\pi$ in view of the ranges where $\widetilde{f}(y), H(y)$ and $L(iy)$ are holomorphic. Note that $\kappa_1\neq0$, otherwise $f(z)$ would continue to an entire function by 1-periodicity, thus contradicting Lemma 6.1.

\smallskip
Now we are ready to conclude the proof. We first note that if $-\pi\rho < \arg(z)<\pi$ then also $-\pi\rho < \arg(z+1), \arg\big(\frac{z}{z+1}\big)<\pi$. Then we apply \eqref{6-18} with $z+1$ in place of $z$, and recalling the 1-periodicity of $f(z)$ we get
\begin{equation}
\label{6-19}
\frac{1}{\kappa_1} \widetilde{\psi}(z+1) = \frac{1}{\kappa_1}f(z) - (z+1)^{-2\mu-1} f\Big(-\frac{1}{z+1}\Big).
\end{equation}
Subtracting \eqref{6-19} from \eqref{6-18}, using the 1-periodicity in the resulting terms $f(-1/(z+1))$ and $f(-1/z)$, and then applying again \eqref{6-18} with $z/(z+1)$ in place of $z$ we finally obtain that
\begin{equation}
\label{6-20}
\begin{split}
 \widetilde{\psi}(z) - \frac{1}{\kappa_1} \widetilde{\psi}(z+1) &= \Big(1-\frac{1}{\kappa_1}\Big) f(z) - \kappa_1 z^{-2\mu-1} f(-1/z) + (z+1)^{-2\mu-1} f\Big(-\frac{1}{z+1}\Big) \\
 &= \Big(1-\frac{1}{\kappa_1}\Big) f(z) + (z+1)^{-2\mu-1} \widetilde{\psi}\Big(\frac{z}{z+1}\Big).
 \end{split}
\end{equation}
Suppose now that $\kappa_1\neq1$. Then \eqref{6-20} gives an expression of $f(z)$ in terms of $ \widetilde{\psi}(z)$, $ \widetilde{\psi}(z+1)$ and $\widetilde{\psi}(z/(z+1))$. But, in particular, $\widetilde{\psi}(z)$ is holomorphic for $|\arg(z)|<\pi\rho$, and clearly $|\arg(z+1)|,|\arg(z/(z+1))|<\pi\rho$ for $z$ in that region. Therefore, such an expression shows that $f(z)$ is holomorphic for $|\arg(z)|<\pi\rho$, thus it is entire by 1-periodicity, a contradiction by Lemma 6.1. Hence 
\begin{equation}
\label{6-21}
\kappa_1=1.
\end{equation}
For future reference we summarize here that thanks to \eqref{6-7},\eqref{6-18} and \eqref{6-21} we have
\begin{equation}
\label{6-22}
\widetilde{\psi}(z)=\psi(z) \quad \text{and} \quad \psi(z) = Q_H(z) + P_H(z) + L(z),
\end{equation}
where
\[
P_H(z) = \frac{(-iz)^{-2\mu-1}}{\Ga(-1/2-\mu)} \widetilde{f}(-iz)
\]
is holomorphic for $|\arg(z)|<\pi$ since we already proved that $\widetilde{f}(y)$ is holomorphic for $y\in\CC\setminus [0,i\infty)$. Thus $\psi(z)$ is holomorphic for $-\pi\rho < \arg(z)<\pi$, and \eqref{6-20} becomes the three-term functional equation
\begin{equation}
 \label{6-23}
\psi(z) = \psi(z+1) + (z+1)^{-2\mu-1} \psi\big(\frac{z}{z+1}\big).
\end{equation}

\smallskip
Finally, we exploit \eqref{6-23} to extend the range where $\psi(z)$ is holomorphic. Indeed, by elementary geometry, we have
\[
\arg(z) < \arg(z+1), \arg\big(\frac{z}{z+1}\big) < 0
\]
for $-\pi < \arg(z) <0$. More precisely, given $R,\ep>0$, if $z\in\C(R,\ep) = \{|z|\leq R:\ -\pi+\ep < \arg(z) < -\ep\}$ then 
\[
\arg(z+1), \arg\big(\frac{z}{z+1}\big) > \arg(z)+\delta(R,\ep)
\]
for some $\delta(R,\ep)>0$. Thus, starting with $z\in\C_\rho(R,\ep) = \{|z|\leq R:\ -\pi\rho < \arg(z) < -\ep\}$, \eqref{6-23} gives step-by-step analytic continuation of $\psi(z)$ to $\C(R,\ep)$, and hence to $|\arg(z)|<\pi$ since $R$ and $\ep$ are arbitrary. The result follows now from \eqref{6-22} since $L(z)$ is holomorphic for $|\arg(z)|<\pi$.  \qed

\medskip
{\bf Proposition 6.2.} {\sl Let $F$ be as in Theorem $1.1$ with $\chi_F\leq 0$.  Then with the notation at the beginning of the section we have
\[
\psi(z) = Q_M(z) + P_M(z) + L(z)
\]
and $Q_M(z)$ is holomorphic for $|\arg(z)|<\pi$.}

\medskip
{\it Proof.} We follow the steps in the proof of Proposition 6.1 till equation \eqref{6-9}. Then, in view of Lemma 5.3, by \eqref{6-4} we have
\begin{equation}
\label{6-24}
\begin{split}
f(iy) &=  \frac{\omega_F}{2\pi i} \int_{(c_1)} (2\pi)^{-s} \Gamma(s) \frac{\ga(1-s+\lambda)}{\ga(s-\lambda)} F(1-s+\lambda) y^{-s} \d s + L(iy) \\
&\hskip1cm + \frac{\omega_F}{2\pi i} \int_{(c_1)} (2\pi)^{-s} \Gamma(s) \big(R(s-\lambda)-1\big) \frac{\ga(1-s+\lambda)}{\ga(s-\lambda)} F(1-s+\lambda) y^{-s} \d s \\
&= \widetilde{P}_M(iy) + L(iy) +  Q_M(iy),
\end{split}
\end{equation}
say. In $\widetilde{P}_M(iy)$ we make the substitution $1-s+\lambda \mapsto s-\lambda$. Thus, recalling \eqref{6-2} and writing 
\begin{equation}
\label{6-25}
T(s) = (2\pi)^{2s-1-2\lambda} \frac{\Gamma(1-s+2\lambda)\ga(s-\lambda)}{\Gamma(s) \ga(1-s+\lambda)}, \quad \omega(\la) =  \omega_F e^{i\frac{\pi}{2}(1+2\lambda)},
\end{equation}
we obtain that
\begin{equation}
\label{6-26}
\begin{split}
\widetilde{P}_M(iy) & = \frac{\omega_F}{2\pi i} \int_{(c_2)} (2\pi)^{s-1-2\lambda} \Gamma(1-s+2\lambda) \frac{\ga(s-\lambda)}{\ga(1-s+\lambda)} F(s-\lambda) y^{s-1-2\lambda} \d s \\
&= \omega(\lambda) \frac{(iy)^{-1-2\lambda}}{2\pi i} \int_{(c_2)} (2\pi)^{-s} \Gamma(s) F(s-\lambda) \big(\frac{1}{y}\big)^{-s} T(s) \d s.
\end{split}
\end{equation}
Now, thanks to \eqref{6-8} and \eqref{6-25}, recalling that $c_2-\Re(\la)>1$ we rewrite \eqref{6-26} in the form
\begin{equation}
\label{6-27}
\begin{split}
\widetilde{P}_M(iy) &= \frac{(iy)^{-1-2\lambda}}{2\pi i} \int_{(c_2)} (2\pi)^{-s} \Gamma(s) F(s-\lambda) \big(\frac{1}{y}\big)^{-s} \d s \\
& \hskip1cm + \frac{(iy)^{-1-2\lambda}}{2\pi i} \int_{(c_2)} (2\pi)^{-s} \Gamma(s) F(s-\lambda) \big(\frac{1}{y}\big)^{-s} \big(\omega(\la)T(s)-1 \big) \d s \\
&=  (iy)^{-2\lambda-1} f\big(-1/(iy)\big) + P_M(iy);
\end{split}
\end{equation}
we use this notation since below we prove that $P_M(iy)$ can be written in the form the form \eqref{6-3}. Hence from \eqref{6-7},\eqref{6-24} and \eqref{6-27} we have by analytic continuation that
\begin{equation}
\label{6-28}
\psi(z) = Q_M(z) + P_M(z) +L(z).
\end{equation}
Thus Proposition 6.2 follows from \eqref{6-28} once we show that $P_M(z)$ can be written as in \eqref{6-3} and $Q_M(z)$ is holomorphic for $|\arg(z)|<\pi$.

\smallskip
To this end we first rewrite $P_M(iy)$ as
\begin{equation}
\label{6-29}
\begin{split}
P_M(iy) =  \frac{(iy)^{-1-2\lambda}}{2\pi i} &\int_{(c_2)} \gamma(s-\lambda) F(s-\lambda) \big(\frac{2\pi}{y}\big)^{-s}\\
&  \hskip-.5cm \times \Big(\omega(\la)(2\pi)^{2s-1-2\lambda}  \frac{\Gamma(1-s+2\lambda)}{\ga(1-s+\lambda)}  - \frac{\Gamma(s)}{\ga(s-\lambda)} \Big) \d s.
\end{split}
\end{equation}
But, in view of \eqref{5-1}, using the duplication and reflection formulae for the $\Ga$ function and the fact that $\ep\in\{0,1\}$  we have
\[
\begin{split}
\frac{\Gamma(1-s+2\lambda)}{\ga(1-s+\lambda)}  &= \frac{\Gamma(1-s+2\lambda)}{\pi^{s-\la-1} \Ga\big(\frac{1-s+\ep+2\la}{2}\big) \Ga\big(\frac{1-s+\ep}{2}\big)} = \pi^{1/2+\la-s} 2^{-s+2\la} \frac{ \Ga\big(\frac{2-s-\ep+2\la}{2}\big)}{\Ga\big(\frac{1-s+\ep}{2}\big)} \\
&= \pi^{1/2+\la-s} 2^{-s+2\la} \frac{ \Ga\big(\frac{2-s-\ep+2\la}{2}\big) \Ga\big(\frac{1+s-\ep}{2}\big)}{\Ga\big(\frac{1-s+\ep}{2}\big) \Ga\big(\frac{1+s-\ep}{2}\big)} \\
&=  \pi^{-1/2+\la-s} 2^{-s+2\la} \Ga\big(\frac{2-s-\ep+2\la}{2}\big) \Ga\big(\frac{1+s-\ep}{2}\big) \sin\big(\pi \big(\frac{1-s+\ep}{2}\big)\big)
\end{split}
\]
and
\[
\begin{split}
\frac{\Gamma(s)}{\ga(s-\lambda)} &= \pi^{s-\la-1/2} 2^{s-1} \frac{\Ga\big(\frac{s}{2}\big) \Ga\big(\frac{s+1}{2}\big)}{\Ga\big(\frac{s+\ep}{2}\big) \Ga\big(\frac{s+\ep-2\la}{2}\big)} =  \pi^{s-\la-1/2} 2^{s-1} \frac{\Ga\big(\frac{s+1-\ep}{2}\big)}{\Ga\big(\frac{s+\ep-2\la}{2}\big)} \\
&= \pi^{s-\la-1/2} 2^{s-1}  \frac{\Ga\big(\frac{s+1-\ep}{2}\big) \Ga\big(\frac{2-s-\ep+2\la}{2}\big)}{\Ga\big(\frac{s+\ep-2\la}{2}\big) \Ga\big(1-\frac{s+\ep-2\la}{2}\big)} \\
&=  \pi^{s-\la-3/2} 2^{s-1} \Ga\big(\frac{1+s-\ep}{2}\big) \Ga\big(\frac{2-s-\ep+2\la}{2}\big) \sin\big(\pi\big(\frac{s+\ep-2\la}{2}\big)\big).
\end{split}
\]
Consequently
\begin{equation}
\label{6-30}
\begin{split}
\omega(\la) (2\pi)^{2s-1-2\lambda}  \frac{\Gamma(1-s+2\lambda)}{\ga(1-s+\lambda)}  &- \frac{\Gamma(s)}{\ga(s-\lambda)}  =  \frac{(2\pi)^s}{2\pi^{\la+3/2}}  \Ga\big(\frac{1+s-\ep}{2}\big) \Ga\big(\frac{2-s-\ep+2\la}{2}\big)  \\
&\times \big(\omega(\la) \sin\big(\pi \big(\frac{1-s+\ep}{2}\big) - \sin\big(\pi\big(\frac{s+\ep-2\la}{2}\big)\big),
\end{split}
\end{equation}
and clearly
\[
\begin{split}
\omega(\la) \sin\big(\pi \big(\frac{1-s+\ep}{2}\big) - \sin\big(\pi\big(\frac{s+\ep-2\la}{2}&\big) = \frac{1}{2i} \Big(-\omega_F e^{i\frac{\pi}{2}(2\la-s+\ep)} - \omega_F e^{i\frac{\pi}{2}(2\la+s - \ep)} \\
& - e^{i\frac{\pi}{2}(s+\ep-2\la)} + e^{-i\frac{\pi}{2}(s+\ep-2\la)}\Big).
\end{split}
\]
Recalling that in this case $\omega_F=(-1)^\ep$, see Proposition 5.1, we have
\[
-\omega_F e^{i\frac{\pi}{2}(2\la-s+\ep)} + e^{-i\frac{\pi}{2}(s+\ep-2\la)} = 0,
\]
and hence
\begin{equation}
\label{6-31}
\omega(\la) \sin\big(\pi \big(\frac{1-s+\ep}{2}\big) - \sin\big(\pi\big(\frac{s+\ep-2\la}{2}\big) = ie^{i\frac{\pi}{2}s} e^{i\frac{\pi}{2}\ep} \cos(\pi\la).
\end{equation}
Gathering \eqref{6-29},\eqref{6-30} and \eqref{6-31} we see that $P_M(iy)$ has the required form \eqref{6-3}.

\smallskip
Finally, from \eqref{6-4} and Lemma 5.3 we have that $Q_M(z)$ is holomorphic for $-\delta < \arg(z) < \pi+\delta$ with some $\delta>0$, hence the same holds for $\psi(z)$ by \eqref{6-28}. But, as in the proof of Proposition 6.1, from the 1-periodicity of $f(z)$ we deduce that $\psi(z)$ satisfies the three-term functional equation \eqref{6-23}, thus the argument at the end of the proof of Proposition 6.1 shows that $\psi(z)$ is holomorphic for $|\arg(z)|<\pi$ also in this case. Proposition 6.2 follows now from \eqref{6-28}. \qed

\smallskip
{\bf Remark 6.1.} The computations in the proof of Proposition 6.2 are similar to those in Lewis-Zagier \cite{Le-Za/2001}, pages 204--205. Apparently there are some slight differences, unimportant in our case, in the final formulae. The function $\psi(z)$ in \eqref{6-7} is a {\it period function} in the sense of Lewis-Zagier \cite{Le-Za/2001}. The fact that $\psi(z)$ in  \eqref{6-28} is holomorphic for $|\arg(z)|<\pi$ follows from the results in Section 4 of Chapter 3 of \cite{Le-Za/2001}; the above argument gives an independent direct proof, which works also for certain functional equations more general than \eqref{6-23}. \qed

\medskip
\section{Conclusion of the proof}

\smallskip
Now we are ready to conclude the proof of Theorem 1.1. Our aim is to show that the function $R(s)$ in Corollary 5.1 is constant, say $R(s)=\eta$. Indeed, if this is the case then $S_F(s) = \eta S_\ga(s)$, hence Corollary 5.1 implies that $F$ satisfies the functions equation
\[
\gamma(s) F(s) = \omega \gamma(1-s) F(1-s) \quad \text{with} \quad \omega = \omega_F \eta,
\]
where $\ga$ is the virtual $\ga$-factor associated with $F$. We are therefore in the situation of Lemma 5.1. More precisely, by \eqref{5-2} and Proposition 5.1, if $\chi_F>0$ then case (i) of Lemma 5.1 holds, otherwise we are in case (ii). Moreover, if $\chi_F=0$ then $\kappa=0$ and hence by Satz 2 of Maass \cite{Maa/1949} we have that $F(s)=\zeta(s)^2$. Theorem 1.1 thus follows. Note, as a side remark, that computing the invariant $\omega_F$ by means of \eqref{2-3}, starting with the above $\ga$-factor and $\omega$-datum, we obtain that $\omega_F=\omega_F\eta$, thus $\eta=1$. 

\smallskip
Suppose, by contradiction, that $R(s)$ is not constant. Then by Lemma 5.3 we have, with the notation in \eqref{5-14}, that $N\geq3$; this will lead to a contradiction.

\smallskip
Suppose first that $\chi_F>0$, i.e. we are in the Hecke case, and let $Q_H(z)$ be as in Proposition 6.1. Applying the reflection formula inside the integral in \eqref{6-5} we rewrite $Q_H(z)$ as
\begin{equation}
\label{7-1}
Q_H(z) = (2\pi)^{-\mu} (-iz)^{-\mu-1} \frac{1}{2\pi i} \int_{(1+\delta)} \frac{S_F(s) + 2\omega_F\cos(\pi s)}{S_F(s) \cos(\pi s)} \sin(\pi(s+\mu)) G_H(z,s) \d s,
\end{equation}
where, in view of \eqref{5-1},
\[
G_H(z,s) = \ga(s) F(s) (-iz)^{-s}.
\]
Next we split the integral in \eqref{7-1} into the sum of the integrals over $(1+\delta -i\infty,1+\delta]$ and $[1+\delta,1+\delta+i\infty)$, thus obtaining that
\begin{equation}
\label{7-2}
Q_H(z) = Q_H^-(z) + Q_H^+(z),
\end{equation}
say. Moreover, standard estimates show that the integrand in $Q_H^+(z)$ is
\[
\ll (t+1)^c e^{(\arg(z)-\pi)t}
\]
with some $c>0$, thus
\begin{equation}
\label{7-3}
Q_H^+(z) \ \text{is holomorphic for} \ |\arg(z)|<\pi.
\end{equation}

\smallskip
Now we observe that by \eqref{5-14}-\eqref{5-16}
\[
S_F(s) + 2\omega_F\cos(\pi s) = \sum_{j=1}^{N-1} a_j e^{i\pi\omega_js},
\]
hence we further split $Q_H^-(z)$ into
\begin{equation}
\label{7-4}
Q_H^-(z) = A_H(z) + B_H(z)
\end{equation}
with
\[
A_H(z) =  (2\pi)^{-\mu} (-iz)^{-\mu-1} \sum_{j=1}^{N-2} \frac{a_j}{2\pi i} \int_{1+\delta-i\infty}^{1+\delta} \frac{e^{i\pi \omega_js}}{S_F(s)\cos(\pi s)} \sin(\pi(s+\mu)) G_H(z,s) \d s
\]
and
\[
B_H(z) =  (2\pi)^{-\mu} (-iz)^{-\mu-1} \frac{a_{N-1}}{2\pi i} \int_{1+\delta-i\infty}^{1+\delta} \frac{e^{i\pi \omega_{N-1}s}}{S_F(s)\cos(\pi s)} \sin(\pi(s+\mu)) G_H(z,s) \d s.
\]
The integrand in $A_H(s)$ is
\[
\ll (|t|+1)^c e^{-|t|(\arg(z) + \pi -\pi\omega_j)}
\]
with some $c>0$, hence, in view of \eqref{5-15}, $A_H(s)$ is holomorphic for $-\pi(\min(1-\omega_{N-2},1)) < \arg(z) < \pi$. Thus, thanks to Proposition 6.1, \eqref{7-2},\eqref{7-3} and \eqref{7-4}, we have that
\begin{equation}
\label{7-5}
B_H(z) \ \text{is holomorphic for} \ -\pi(\min(1-\omega_{N-2},1)) < \arg(z) < \pi.
\end{equation}

\smallskip
The last step is to suitably transform the function $B_H(z)$. To this end we first note that as $t\to-\infty$
\begin{equation}
\label{7-6}
\frac{\sin(\pi(s+\mu))}{S_F(s)\cos(\pi s)} =c_1 e^{-i\pi s} \big( 1 + O(e^{-\pi(1-\omega_{N-1})|t|})\big)
\end{equation}
with a certain constant $c_1\neq0$. Therefore
\begin{equation}
\label{7-7}
B_H(z) = c_2 (-iz)^{-\mu-1} \frac{1}{2\pi i} \int_{1+\delta-i\infty}^{1+\delta} e^{i\pi s(\omega_{N-1}-1)} G_H(z,s) \d s + C_H(z)
\end{equation}
with a certain constant $c_2\neq0$, where
\begin{equation}
\label{7-8}
C_H(z) \ \text{is holomorphic for} \ -\pi(\min(2-2\omega_{N-1},1)<\arg(z)<\pi
\end{equation}
thanks to \eqref{7-6}. Next we extend the integral in \eqref{7-7} to the whole line $\si=1+\delta$. Since as $t\to+\infty$ the integrand in \eqref{7-7} is
\[
\ll (|t|+1)^c e^{|t|(\arg(z) -\pi\omega_{N-1})}
\]
with some $c>0$, the corresponding integral is holomorphic for $-\pi <\arg(z) < \pi\omega_{N-1}$. Therefore, writing
\[
I_H(z) =  \frac{1}{2\pi i} \int_{(1+\delta)} e^{i\pi s(\omega_{N-1}-1)} G_H(z,s) \d s,
\]
from \eqref{7-5},\eqref{7-7} and \eqref{7-8} we deduce that $I_H(z)$ is holomorphic for
\begin{equation}
\label{7-9}
-\pi\min(1-\omega_{N-2},2-2\omega_{N-1},1) < \arg(z) < \pi\omega_{N-1}.
\end{equation}
But in view of \eqref{6-8} we have
\[
I_H(z) = \left(-2\pi ie^{i\pi(1-\omega_{N-1})}z\right)^\mu f\left(e^{i\pi(1-\omega_{N-1})}z\right),
\]
hence thanks to \eqref{7-9} and recalling that $\omega_{N-2}<\omega_{N-1}$ we finally conclude that $f(z)$ is holomorphic for $-\rho\pi < \arg(z)<\pi$ with some $\rho>0$. This is a contradiction by Lemma 6.1 thanks to the 1-periodicity of $f(z)$, thus finishing the proof of the Hecke case.

\smallskip
Suppose now that $\chi_F\leq0$, i.e. we are in the Maass case. The proof is similar to the Hecke case, so we will be more sketchy. Recalling that $\la=i\kappa$, let
\begin{equation}
\label{7-10}
G_M(z,s) = \Gamma(s) \frac{\ga(1-s+\lambda)}{\ga(s-\lambda)} F(1-s+\lambda) (-2\pi iz)^{-s}.
\end{equation}
Thanks to Lemma 5.3, \eqref{5-14} and \eqref{5-15}, we rewrite the function $Q_M(z)$ in \eqref{6-4} as
\begin{equation}
\label{7-11}
\begin{split}
Q_M(z)& = \frac{1}{2\pi i} \Big(\int_{c_1-i\infty}^{c_1} + \int_{c_1}^{c_1+i\infty} \Big) \big(R(s-\lambda)-1) G_M(z,s) \d s = Q_M^-(z) + Q_M^+(z) \\
&= Q_M^+(z) + \sum_{j=0}^{N-2} \frac{a_j}{2\pi i}  \int_{c_1-i\infty}^{c_1} \frac{e^{i\pi\omega_j(s-\la)}}{S_\ga(s-\la)} G_M(z,s) \d s \\
& \hskip1cm+ \frac{a_{N-1}}{2\pi i}  \int_{c_1-i\infty}^{c_1} \frac{e^{i\pi\omega_{N-1}(s-\la)}}{S_\ga(s-\la)} G_M(z,s) \d s \\
& \hskip1cm+ \frac{1}{2\pi i}  \int_{c_1-i\infty}^{c_1} \Big(\frac{a_Ne^{i\pi(s-\la)}}{S_\ga(s-\la)}-1\Big) G_M(z,s) \d s \\
&=  Q_M^+(z) + A_M(z) + B_M(z) + C_M(z),
\end{split}
\end{equation}
say. Standard estimates involving Lemma 5.3 show that
\begin{equation}
\label{7-12}
\text{$Q_M^+(z)$ is holomorphic for $|\arg(z)|<\pi$.}
\end{equation}
For $0\leq j \leq N-2$ and $t\to-\infty$ we have
\[
\frac{e^{i\pi\omega_j(s-\la)}}{S_\ga(s-\la)} G_M(z,s) \ll (|t|+1)^c e^{-|t|(\arg(z)+\pi(1-\omega_j))}
\]
for some $c>0$, hence
\begin{equation}
\label{7-13}
\text{$A_M(z)$ is holomorphic for $-\pi\min(1-\omega_{N-2},1)<\arg(z)<\pi$.}
\end{equation}
Moreover, recalling \eqref{5-16},\eqref{5-23} and Proposition 5.1, we observe that the main terms of $a_Ne^{i\pi(s-\la)}$ and $S_\ga(s-\la)$ are equal, hence as $t\to-\infty$
\[
\frac{a_Ne^{i\pi(s-\la)}}{S_\ga(s-\la)}-1 \ll e^{-\pi|t|}.
\]
Therefore
\begin{equation}
\label{7-14}
\text{$C_M(z)$ is holomorphic for $|\arg(z)|<\pi$.}
\end{equation}
Finally, to treat $B_M(z)$ we observe from \eqref{5-23} that as $t\to-\infty$
\[
\frac{1}{S_\ga(s-\la)} = \alpha_1  e^{-i\pi s} + O(e^{-2\pi|t|})
\]
with a certain constant $\alpha_1\neq0$, hence
\begin{equation}
\label{7-15}
B_M(z) = \frac{\alpha_2}{2\pi i} \int_{c_1-i\infty}^{c_1} e^{i\pi s(\omega_{N-1}-1)} G_M(z,s) \d s + D_M(z)
\end{equation}
 with a certain constant $\alpha_2\neq0$ and a function $D_M(z)$ holomorphic for $|\arg(z)|<\pi$. But as $t\to+\infty$ the integrand in \eqref{7-15} is
 \[
 \ll (|t|+1)^c e^{-|t|(\pi\omega_{N-1} -\arg(z))}
 \]
 for some $c>0$. Thus we have that
\begin{equation}
\label{7-16}
 \frac{\alpha_2}{2\pi i} \int_{c_1}^{c_1+i\infty} e^{i\pi s(\omega_{N-1}-1)} G_M(z,s) \d s \quad \text{is holomorphic for $-\pi <\arg(z) < \pi\omega_{N-1}$.}
 \end{equation}
Therefore, gathering \eqref{7-10}-\eqref{7-16}, from Proposition 6.2 we deduce that the integral
\begin{equation}
\label{7-17}
I_M(z)= \frac{1}{2\pi i} \int_{(c_1)} \Gamma(s) \frac{\ga(1-s+\lambda)}{\ga(s-\lambda)} F(1-s+\lambda) (-2\pi i e^{i\pi (1-\omega_{N-1})}  z)^{-s}  \d s
\end{equation}
represents a function holomorphic for $-\pi\min(1-\omega_{N-2},1)<\arg(z)<\pi\omega_{N-1}$.

\smallskip
Thanks to \eqref{5-7}, Proposition 5.1, \eqref{2-10},\eqref{2-1} and Lemma 4.1, recalling that $\ga_F(s)$ denotes the $\ga$-factor of $F$ we have
\[
\begin{split}
\frac{\ga(1-s+\lambda)}{\ga(s-\lambda)} &= \frac{\ga_F(1-s+\lambda)}{\ga_F(s-\lambda)} \frac{S_\ga(s-\la)}{S_F(s-\la)} \\
&=   \frac{\ga_F(1-s+\lambda)}{\ga_F(s-\lambda)}  +  \frac{\ga_F(1-s+\lambda)}{\ga_F(s-\lambda)} \Big(\frac{S_\ga(s-\la)}{S_F(s-\la)} -1 \Big) \\
&=    \frac{\ga_F(1-s+\lambda)}{\ga_F(s-\lambda)}  + O(e^{-\delta|t|})
\end{split}
\]
with some $\delta>0$ by Lemma 5.3. Therefore, replacing $\frac{\ga(1-s+\lambda)}{\ga(s-\lambda)}$ by $\frac{\ga_F(1-s+\lambda)}{\ga_F(s-\lambda)}$ in \eqref{7-17} adds an extra term which is holomorphic at least for $-\pi\min(1-\omega_{N-1}+\delta,1)<\arg(z)<\pi\omega_{N-1}$. Thus, writing 
\[
w=w(z) := -i e^{i\pi (1-\omega_{N-1})}  z, 
\]
the integral
\[
J_M(z) = \frac{1}{2\pi i} \int_{(c_1)} \Gamma(s) \frac{\ga_F(1-s+\lambda)}{\ga_F(s-\lambda)} F(1-s+\lambda) (2\pi w)^{-s}  \d s
\]
represents a function holomorphic for 
\begin{equation}
\label{7-18}
-\pi\min(1-\omega_{N-2},1-\omega_{N-1}+\delta,1)<\arg(z)<\pi\omega_{N-1}.
 \end{equation}
But by \eqref{2-1}
\[
J_M(z) = \frac{\omega_F}{2\pi i} \int_{(c_1)} \Gamma(s) F(s-\lambda) (2\pi w)^{-s}  \d s.
\]
Therefore, recalling the choice of $c_1$ in \eqref{6-2} and shifting the integration line to $\si=c>1+\Re(\la)$, by Mellin's transform we get
\[
\begin{split}
J_M(z) &= \frac{\omega_F}{2\pi i} \int_{(c)} \Gamma(s) F(s-\lambda) (2\pi w)^{-s}  \d s - \omega_F \res_{s=1+\la} \Gamma(s) F(s-\lambda) (2\pi w)^{-s} \\
&= \omega_F f\big(e^{i\pi (1-\omega_{N-1})}  z\big)  - \omega_F \res_{s=1+\la} \Gamma(s) F(s-\lambda) (2\pi w)^{-s}.
\end{split}
\]
Thus $f\big(e^{i\pi (1-\omega_{N-1})}  z\big)$ is holomorphic in the region \eqref{7-18}, hence $f(z)$ is holomorphic for $-\rho\pi<\arg(z)<\pi$ for some $\rho>0$, a contradiction as in the Hecke case. Theorem 1.1 is now proved. \qed

\newpage

\ifx\undefined\bysame{poly}.
\newcommand{\bysame}{\leavevmode\hbox to3em{\hrulefill}\ ,}
\fi

\bigskip
\bigskip
\noindent
Jerzy Kaczorowski, Faculty of Mathematics and Computer Science, A.Mickiewicz University, 61-614 Pozna\'n, Poland and Institute of Mathematics of the Polish Academy of Sciences, 00-956 Warsaw, Poland. e-mail: \url{kjerzy@amu.edu.pl}

\medskip
\noindent
Alberto Perelli, Dipartimento di Matematica, Universit\`a di Genova, via Dodecaneso 35, 16146 Genova, Italy. e-mail: \url{perelli@dima.unige.it}

\end{document}